\documentclass[bj]{imsart}

\RequirePackage[OT1]{fontenc}
\RequirePackage{amsthm,amsmath}
\RequirePackage[colorlinks,citecolor=blue,urlcolor=blue]{hyperref}

\startlocaldefs
\numberwithin{equation}{section}
\theoremstyle{plain}

\theoremstyle{remark}

\endlocaldefs

\usepackage[square,sort,comma,numbers]{natbib}
\usepackage{latexsym}
\usepackage{amssymb}
\usepackage{amsmath}
\usepackage{amsthm}
\usepackage[mathscr]{eucal}
\usepackage{graphicx}
\usepackage{caption}
\usepackage{subcaption}

\newtheorem{Lemma}{Lemma}

\newtheorem{Theorem}[Lemma]{Theorem}

\renewcommand{\qed}{\hfill{\ \ \rule{2mm}{2mm}} \vspace{0.2in}}

\newcommand{\ind}{1\hspace{-2.3mm}{1}}

\begin{document}

\begin{frontmatter}
\title{Minimum spanning trees of random geometric graphs with location dependent weights}
\runtitle{Minimum spanning trees }

\begin{aug}
\author{\fnms{Ghurumuruhan} \snm{Ganesan}\thanksref{a,e1}\ead[label=e1,mark]{gganesan82@gmail.com}}%
%\ead[label=u1,url]{www.foo.com}}

\address[a]{Institute of Mathematical Sciences, Chennai.
\printead{e1}}

\runauthor{G. Ganesan}

\affiliation{Institute of Mathematical Sciences, Chennai}

\end{aug}

%\title{}
%\author{ \textbf{Ghurumuruhan Ganesan}
%\thanks{E-Mail: \texttt{gganesan82@gmail.com} } \\
%EndAName
%\ \\
%Institute of Mathematical Sciences, Chennai.}
\date{}
\maketitle

\begin{abstract}
Consider~\(n\) nodes~\(\{X_i\}_{1 \leq i \leq n}\)
independently distributed in the unit square~\(S,\)
each according to a distribution~\(f.\) Nodes~\(X_i\) and~\(X_j\)
are joined by an edge if the Euclidean distance~\(d(X_i,X_j)\) is less than~\(r_n,\)
the adjacency distance and the resulting random graph~\(G_n\)
is called a random geometric graph~(RGG). We now
assign a location dependent weight to each edge of~\(G_n\)
and define~\(MST_n\) to be the sum of the weights of the minimum spanning trees of
all components of~\(G_n.\) For values of~\(r_n\) above the connectivity regime,
we obtain upper and lower bound deviation estimates for~\(MST_n\)
and~\(L^2-\)convergence of~\(MST_n\) appropriately
scaled and centred. %converges to zero in mean as~\(n \rightarrow \infty.\)

%\(\frac{1}{\sqrt{n}}(MST_n - \mathbb{E}MST_n)\)

 %and analogous results hold for minimum spanning trees and paths.

%In this paper, we study the structure of left-right crossings
%of the random geometric graph \(G = G(n,r_n)\) of \(n\) nodes
%uniformly distributed in \(S = [0,1]^2\) with \(r_n = \epsilon\sqrt{\frac{\log{n}}{n}}\)
%for some \(\epsilon > 0.\) Tiling \(S\) horizontally and
%vertically into rectangles of length \(1\) and width \(Mr_n,\) we
%show that each rectangle has a left-right crossing of edges with
%high probability if \(M\) is sufficiently large.
%We call the resulting subgraph to be a ``backbone" of \(G.\)

%The techniques we use to construct the backbone has quite a few applications.
%As a first, we show that the diameter of second largest component in \(G\)
%is \(O(1)\) with high probability. Secondly,
%\vspace{0.1in} \noindent \textbf{Key words:} .

\vspace{0.1in} \noindent \textbf{AMS 2000 Subject Classification:} Primary:
60J10, 60K35; Secondary: 60C05, 62E10, 90B15, 91D30.
\end{abstract}


\end{frontmatter}

\bigskip

\setcounter{equation}{0}
\renewcommand\theequation{\thesection.\arabic{equation}}
\section{Introduction}\label{intro}
The study of minimum spanning trees (MSTs) of a graph arise in many applications
and many analytical results have been derived regarding the weight of the MST
for various types of weighted graphs. For MSTs of complete Euclidean graphs
with edge weights being Euclidean length,~\cite{beard} study convergence of the weight
of the travelling salesman path and describe modifications
that allow for the study of MSTs appropriately scaled and centred.
Also,~\cite{pen}~\cite{steele2}~\cite{yuk}
study MSTs of complete Euclidean graphs where
the weights of the edges grow as a power of the Euclidean length
of the edges. Asymptotic convergence of the MSTs (in probability
and in the sense of CLT) together with auxiliary results
are studied as the number of nodes~\(n \rightarrow \infty\)
by estimating the expected weight using the bounded degree property
of MSTs (see also~\cite{steele}~\cite{yuk2}~\cite{kest} and references therein for
more details). For material on the algorithmic and application oriented aspects of MSTs, we refer
to~\cite{djau}~\cite{jothi} and references therein.

In this paper, we consider random geometric graphs (RGGs)~\cite{pen2}~\cite{pen3}
where nodes are distributed randomly across the unit square and nodes
close enough to each other are connected to each other by edges.
Also, the weight of an edge in the RGG might depend on the \emph{individual} locations
of the endvertices.

The scenario described above arises frequently in wireless networks
and for example, suppose it is required to establish a fully connected wireless communication network among the set of nodes distributed
randomly in a certain geographical area. The cost of setting up a communication link between any two nodes (e.g. length of the cables etc) is usually assumed to be directly proportional to the distance between the nodes~\cite{gupta}. Additionally, the cost could also depend on the location of the nodes since it might happen that some parts of area are ``remote" and so it might cost more to install links among nodes located in these parts.

In the rest of this section, we describe the model under consideration and state our main results regarding the MSTs of RGGs with location dependent weights. In our main result (Theorem~\ref{mst_thm}), we obtain deviation and variance estimates for the weight of the MST, appropriately scaled and centred.

%Using our proof techniques, we also obtain variance estimates

\subsection*{Model Description}
Let~\(n\) nodes~\(\{X_i\}_{1 \leq i \leq n}\) be independently
distributed in the unit square~\(S,\) each according to a density~\(f\) satisfying
\begin{equation}\label{f_eq}
\epsilon_1 \leq \inf_{x \in S} f(x) \leq \sup_{x \in S} f(x) \leq \epsilon_2
\end{equation}
for some positive constants~\(\epsilon_1 \leq 1 \leq \epsilon_2,\) since~\(\int_{S} f(x)dx = 1.\)
Throughout constants do not depend on~\(n.\) The Euclidean distance between nodes~\(X_i\) and~\(X_j\)
is denoted by~\(d(X_i,X_j)\) and nodes~\(X_i\) and~\(X_j\)
are connected by an edge if~\(d(X_i,X_j) < r_n.\) The term~\(r_n\)
is called the adjacency distance and the resulting random graph~\(G_n\)
is called a random geometric graph (RGG) (Penrose~(2003)).

Let~\(Y_1,\ldots,Y_t \subset \{X_k\}_{1 \leq k \leq n}\) be~\(t \geq 2\) distinct nodes. A path~\({\cal Q} = (Y_1,\ldots,Y_t)\) is a subgraph of~\(G_n\)
with vertex set~\(\{Y_{j}\}_{1 \leq j \leq t}\)
and edge set~\(\{(Y_{j},Y_{{j+1}})\}_{1 \leq j \leq t-1}.\)
The nodes~\(Y_1\) and~\(Y_t\) are said to be \emph{connected} by
edges of the path~\({\cal Q}.\) For~\(t \geq 3,\) the subgraph~\({\cal C}  = (Y_1,Y_2,\ldots,Y_t,Y_1)\) with
vertex set~\(\{Y_{j}\}_{1 \leq j \leq t}\)
and edge set~\(\{(Y_{j},Y_{{j+1}})\}_{1 \leq j \leq t-1} \cup \{(Y_t,Y_1)\}\)
is said to be a \emph{cycle}.

%For an edge~\(e \in K_n\) let~\(l(e)\) denote the length of~\(e.\)

A subgraph~\({\cal T}\) of~\(G_n\)
with vertex set~\(\{Y_i\}_{1 \leq i \leq t}\) and edge set~\(E_{\cal T}\)
is said to be a \emph{tree} if
the following two conditions hold:\\
\((1)\) The graph~\({\cal T}\) is connected; i.e., any two nodes
in~\({\cal T}\) are connected by a path containing only edges in~\(E_{\cal T}.\)\\
\((2)\) The graph~\({\cal T}\) is acyclic; i.e., no subgraph of~\({\cal T}\) is a cycle.\\
The tree~\({\cal T}\) is said to be a \emph{spanning tree} of a component~\({\cal C}\) of the graph~\(G_n\) if~\({\cal T}\)
contains all the nodes of~\({\cal C}.\) If~\({\cal T}\) contains all the~\(n\) nodes~\(\{X_k\}_{1 \leq k \leq n}\)
then~\({\cal T}\) is said to be a spanning tree of the RGG~\(G_n.\)

We are interested in studying minimum spanning trees of~\(G_n\) each of whose edges
is assigned a random weight as described below.

\subsubsection*{Edge weights}
Let~\(\xi : S \times S \rightarrow (0,\infty)\) be any measurable function satisfying~\(\xi(x,y) = \xi(y,x)\) for all~\(x,y \in S.\) For~\(1 \leq i < j \leq n\) we define~\[w(X_i,X_j) := d^{\alpha}(X_i,X_j)\cdot \xi(X_i,X_j)\] to be the weight of the edge~\((X_i,X_j),\) associated with the (deterministic) edge weight factor~\(\xi\) and exponent~\(\alpha > 0.\) Throughout we assume that~\(\xi(x,y) \in [\xi_{min}, \xi_{max}]\) for some positive constants~\(\xi_{min}\) and~\(\xi_{max}.\) For a component~\({\cal C}\) of the RGG~\(G_n\) with vertex set~\(\{Y_1,\ldots, Y_t\}\) and for a spanning tree~\({\cal T}\) of~\({\cal C},\) the \emph{weight} of~\({\cal T}\) is the sum of the weights of the edges in~\({\cal T};\) i.e.,
\begin{equation}\label{len_cyc_def}
W({\cal T})  = W({\cal T},\xi) := \sum_{e \in {\cal T}} w(e) = \frac{1}{2} \sum_{i=1}^{t} w(Y_i,{\cal T}),
\end{equation}
where~\(w(Y_i,{\cal T})\) is the sum of the weights of edges in~\({\cal T}\) containing~\(Y_i\) as an endvertex. Let
\begin{equation}\label{min_weight_tree}
MST_n = MST_n\left(G_n,\xi\right) := \sum_{{\cal C} \in G_n}\min_{{\cal T}} W({\cal T})
\end{equation}
where the sum is taken over all components of~\(G_n\) and for each component~\({\cal C},\) the minimum is taken over all spanning trees~\({\cal T}\) of~\({\cal C}.\) If~\(G_n\) is connected, then we denote the spanning tree~\({\cal T}_n\) with weight~\(MST_n\) to be the \emph{minimal spanning tree} (MST). If there is more than one choice for~\({\cal T}_n,\) we choose one according to a deterministic rule.

We emphasize here that the edge weight factor~\(\xi\) is deterministic and the only randomness in~\(MST_n\) comes from the location of the nodes~\(\{X_i\}_{1 \leq i \leq n}.\) Let~\(\epsilon_1,\epsilon_2\) be as in~(\ref{f_eq}) and set~\(\delta = \delta(\alpha) = \epsilon_1\) if the edge weight exponent~\(\alpha \leq 1\) and~\(\delta = \epsilon_2\) if~\(\alpha > 1.\) For~\(A > 0\) we define~\(C_1(A)=C_1(A,\epsilon_1,\epsilon_2,\alpha)\)
and~\(C_2(A) = C_2(A,\epsilon_1,\epsilon_2,\alpha)\) as
\begin{eqnarray}\label{c12def}
C_1(A) &:=& \frac{1}{2}\xi_{min}A^{\alpha-2}(1-e^{-\epsilon_1 A^2})e^{-8\epsilon_2 A^2} \text{ and } \nonumber\\
C_2(A)  &:=& \xi_{max}(2A)^{\alpha}\left(1 + \frac{\mathbb{E}\tilde{T}^{\alpha}}{A^2}\right),
\end{eqnarray}
where~\(\tilde{T}\) is a geometric random variable with success parameter~\(p = 1-e^{-\delta A^2}.\) We have the following result.
\begin{Theorem}\label{mst_thm} Suppose the following conditions are satisfied:\\
\((i)\) The adjacency distance~\(r_n \geq \sqrt{\frac{M\log{n}}{n}}\) and~\(r_n \longrightarrow 0\) as~\(n \rightarrow \infty,\) where\\\(M > \left(\frac{1600}{\epsilon_1}\right)\) is a constant.\\
\((ii)\) The edge weight exponent~\(\alpha > 0\) and the edge weight factor~\(\xi(x,y) \in [\xi_{min}, \xi_{max}]\) for all~\(x,y\) and some positive constants~\(\xi_{min},\xi_{max}.\)\\
There is a constant~\(C> 0\) such that
\begin{equation}\label{var_mst_est_main}
var\left(\frac{MST_{n}}{n^{1-\frac{\alpha}{2}}}\right) \leq Cr_n^2(nr_n^2)^{\alpha}.
\end{equation}
Moreover, for every~\(A > 0\) there exists~\(A_n \in \left[A+ \frac{1}{(\log{n})^{1/4}}, A + \frac{2}{(\log{n})^{1/4}}\right)\) and a constant~\(D > 0\) such that
\begin{equation}\label{mst_low_bounds}
\mathbb{P}\left(MST_n \geq C_1(A_n) n^{1-\frac{\alpha}{2}}\left(1-\frac{36\sqrt{A_n}}{n^{1/4}}\right)\right) \geq 1-\frac{1}{n^2}
\end{equation}
and
\begin{equation}\label{mst_up_bounds}
\mathbb{P}\left(MST_n \leq C_2(A_n) n^{1-\frac{\alpha}{2}}\left(1+ \frac{1}{n^{1/17}}\right)\right) \geq 1-\frac{D}{n^{5/4}}
\end{equation}
for all~\(n\) large. Consequently,
\begin{equation}\label{exp_mst_bound}
C_1(A_n) \left(1-\frac{37\sqrt{A_n}}{n^{1/4}}\right)\leq \mathbb{E}\left(\frac{MST_n}{n^{1-\frac{\alpha}{2}}}\right) \leq C_2(A_n)\left(1+ \frac{2}{n^{1/17}}\right)
\end{equation}
for all~\(n\) large.
\end{Theorem}
We remark here that we have chosen~\(r_n\) to be at least of the order of~\(\sqrt{\frac{\log{n}}{n}}\) so that the resulting RGG is connected with high probability (see discussion following~(\ref{z_tot_def}) in Section~\ref{pf_mst_var}). From the variance estimate~(\ref{var_mst_est_main}) in Theorem~\ref{mst_thm}, we see that if
\begin{equation}\label{rn_cond2}
n^{\frac{1}{2} \cdot \frac{\alpha}{1+\alpha}}\cdot r_n \longrightarrow 0
\end{equation}
as~\(n \rightarrow \infty,\) then~\[\frac{MST_n-\mathbb{E}MST_n}{n^{1-\frac{\alpha}{2}}} \longrightarrow 0\] in~\(L^2\) as~\(n \rightarrow \infty.\) For example if~\(\alpha = 1\) and~\(r_n = \frac{1}{n^{1/3}},\) then both condition~\((i)\) in Theorem~\ref{mst_thm} and~(\ref{rn_cond2}) are satisfied.

We use~(\ref{mst_up_bounds}) to evaluate the asymptotic numeric bounds for~\(\mathbb{E}MST_n\) as follows. Recalling the definition of~\(A_n\) in the statement of Theorem~\ref{mst_thm}, we prove in the Appendix that for~\(A > 0,\) the limits
\begin{equation}\label{c_bounds}
\lim_{n} C_i(A_n) = C_i(A) \text{ for } i= 1,2.
\end{equation}
Therefore defining
\[\beta_{up} := \inf_{A > 0} C_2(A) \text{ and } \beta_{low} := \sup_{A > 0} C_1(A),\] we get from~(\ref{exp_mst_bound}) that
\[\beta_{low} \leq \liminf_{n} \mathbb{E}\left(\frac{MST_n}{n^{1-\frac{\alpha}{2}}}\right)   \leq \limsup_{n} \mathbb{E}\left(\frac{MST_n}{n^{1-\frac{\alpha}{2}}}\right)  \leq \beta_{up}.\] For example, for the homogenous case of~\(\epsilon_1 = \epsilon_2 = \xi_{min} = \xi_{max} = 1\) and~\(\alpha = 1,\) we have that
\[\beta_{low} \approx 0.0735633 \text{ and } \beta_{up} \approx 4.46256.\]

We briefly outline the methods in the proof of Theorem~\ref{mst_thm}. The variance bound~(\ref{var_mst_est_main}) is obtained via the martingale difference method together with one node difference estimates that measures the change in MST lengths after adding or removing a single node, retaining the same adjacency distance.

To prove the deviation estimates, we use Poissonization and prove results for the Poissonized RGG and then dePoissonize to obtain the corresponding results for the Binomial RGG. In the rest of the two paragraphs we simply denote the Poissonized RGG as RGG. First we show that the event~\(E_{conn}\) that the (Poissonized) RGG is connected occurs with high probability, i.e., with probability converging to one as~\(n \rightarrow \infty\) and consider only subevents of~\(E_{conn}.\) For the lower deviation estimate, we tile the unit square into roughly~\(n\) small squares~\(\{R_k\}\) of side length of order~\(\frac{1}{\sqrt{n}}\) each and look for special type of occupied squares, whose neighbouring squares are all vacant. This results in a vacant annulus around such squares and since the RGG is connected, nodes within these special squares must have edges crossing over the vacant annulus to reach other nodes. This allows us to construct edges of length of order~\(\frac{1}{\sqrt{n}}\) in the MST and determining that there are order of~\(n\) such special squares with high probability, we obtain the lower deviation bound.

For the upper deviation bound, we join nodes within each square in~\(\{R_k\}\) to get a collection of subtrees. We then join all these subtrees together by adding extra edges to get an spanning tree whose weight is no more than the MST. Further the weight of this constructed tree has two parts: (1) the sum length of the edges of the subtrees, each of which is of order~\(\frac{1}{\sqrt{n}}\) by construction and \((2)\) the sum length of the extra edges. To estimate this second sum, we use coupling with homogenous Poisson process and deviation estimates for Geometric random variables to obtain the desired upper bound. For more details, we refer to Section~\ref{pf_mst_thm_dev} and finally, we remark that the methods used in this paper can also be used to analyze RGGs in~\(d-\)dimensions and in regular shapes other than the unit square.

% ns_n^{2\alpha}\xi_{max}^2 \left((4\pi^2\epsilon_2^2 + 2\pi\epsilon_2)\omega^2_n + (600\epsilon_2 ns_n^2)^2\right),

%we get that if the adjacency distance~\(r_n\) satisfies~\(r_n \cdot n^{\frac{\alpha}{2+2\alpha}} \longrightarrow 0\) as~\(n \rightarrow \infty\) then~\(\frac{1}{n^{1-\frac{\alpha}{2}}}(MST_n - \mathbb{E}MST_n) \longrightarrow 0\) in probability as~\(n \rightarrow \infty.\) One such example is the sequence~\(r_n = \frac{1}{n^{\theta}}\) where~\(\frac{\alpha}{2+2\alpha} < \theta < \frac{1}{2}.\) For such specific values of~\(r_n,\) however, we have stronger a.s.\ convergence as described in the following result.

The paper is organized as follows. In Section~\ref{pf_mst_var}, we prove the variance estimate~(\ref{var_mst_est_main}) in Theorem~\ref{mst_thm} and  in Section~\ref{pf_mst_thm_dev}, we prove the deviation estimates~(\ref{mst_low_bounds}) and~(\ref{mst_up_bounds}) and the expectation bounds~(\ref{exp_mst_bound}) in Theorem~\ref{mst_thm}. %In Section~\ref{pf_mst_conv}, we prove Theorem~\ref{mst_thm_conv}.

%CHK REF ETC.. AND ALL +eTC...

%WRT SMALL PROF OF COR +etC...

%In the end, we convert the estimates to the Binomial process.

%Arguing iteratively as in the proof of~\((h2),\) we get
%\begin{eqnarray}
%\left|MST_{k} - MST_{n^2}\right| &\leq& C_1r_{n^2} (\log{n^2})(k-n^2)\ind\left(Y_{tot}(n^2)\right) \nonumber\\
%&&\;\;\;\;\;\;+\;\;\;(k-n^2) n^2\sqrt{2} \ind\left(Y_{tot}^c(n^2)\right)\label{d_est_mst}
%\end{eqnarray}
%GOOTHALS HERE...

\setcounter{equation}{0}
\renewcommand\theequation{\thesection.\arabic{equation}}
\section{Proof of the variance estimate in Theorem~\ref{mst_thm}}\label{pf_mst_var}
Throughout we use the following standard deviation estimates~\cite{alon}
for sums of independent Poisson and Bernoulli random variables. Suppose~\(W_i, 1 \leq i \leq m\) are independent Bernoulli random variables satisfying~\(\mu_1 \leq \mathbb{P}(W_1=1) = 1-\mathbb{P}(W_1~=~0) \leq \mu_2.\) For any~\(0 < \epsilon < \frac{1}{2},\)
\begin{equation}\label{std_dev_up}
\mathbb{P}\left(\sum_{i=1}^{m} W_i > \mu_2(1+\epsilon) \right) \leq \exp\left(-\frac{\epsilon^2}{4}\mu_2\right)
\end{equation}
and
\begin{equation}\label{std_dev_down}
\mathbb{P}\left(\sum_{i=1}^{m} W_i < \mu_1(1-\epsilon) \right) \leq \exp\left(-\frac{\epsilon^2}{4}\mu_1\right)
\end{equation}
Estimates~(\ref{std_dev_up}) and~(\ref{std_dev_down}) also hold if~\(\{W_i\}\) are independent Poisson random variables with~\(\mu_1 \leq \mathbb{E}W_1 \leq \mu_2.\)

\subsection*{One node difference estimates}
In this subsection we find estimates for changes in the length of the MST upon adding or removing a single node. By definition, the adjacency distance~\(r_n\) depends on the number of nodes~\(n\) in consideration and so strictly speaking, adding or removing a node changes the adjacency distance and hence the underlying RGG. However, for our purposes we obtain below difference estimates \emph{retaining} the same adjacency distance.

As in~(\ref{min_weight_tree}), let~\(MST_{n+1}\) be the sum of the weights of the MSTs of all the components of the RGG~\(G_{n+1},\) formed by the nodes~\(\{X_i\}_{1 \leq i \leq n+1}\) with adjacency distance~\(r_{n+1}.\) For~\(1 \leq i \leq n+1,\) let~\(G(i,r_{n+1})\) be the RGG formed by the~\(n\) nodes~\(\{X_j\}_{1 \leq j \neq i \leq n+1}\) with the \emph{same} adjacency distance~\(r_{n+1}\) and let~\(MST_n(i,r_{n+1})\) be the sum of the weights of the MSTs of all the components of the RGG~\(G(i,r_{n+1}),\) as defined in~(\ref{min_weight_tree}). We are interested in estimating~\(|MST_{n+1} - MST(i,r_{n+1})|\) for any~\(1 \leq i \leq n+1.\)

\begin{figure}[tbp]
\centering
%\fbox{
\includegraphics[width=3in, trim= 20 200 50 110, clip=true]{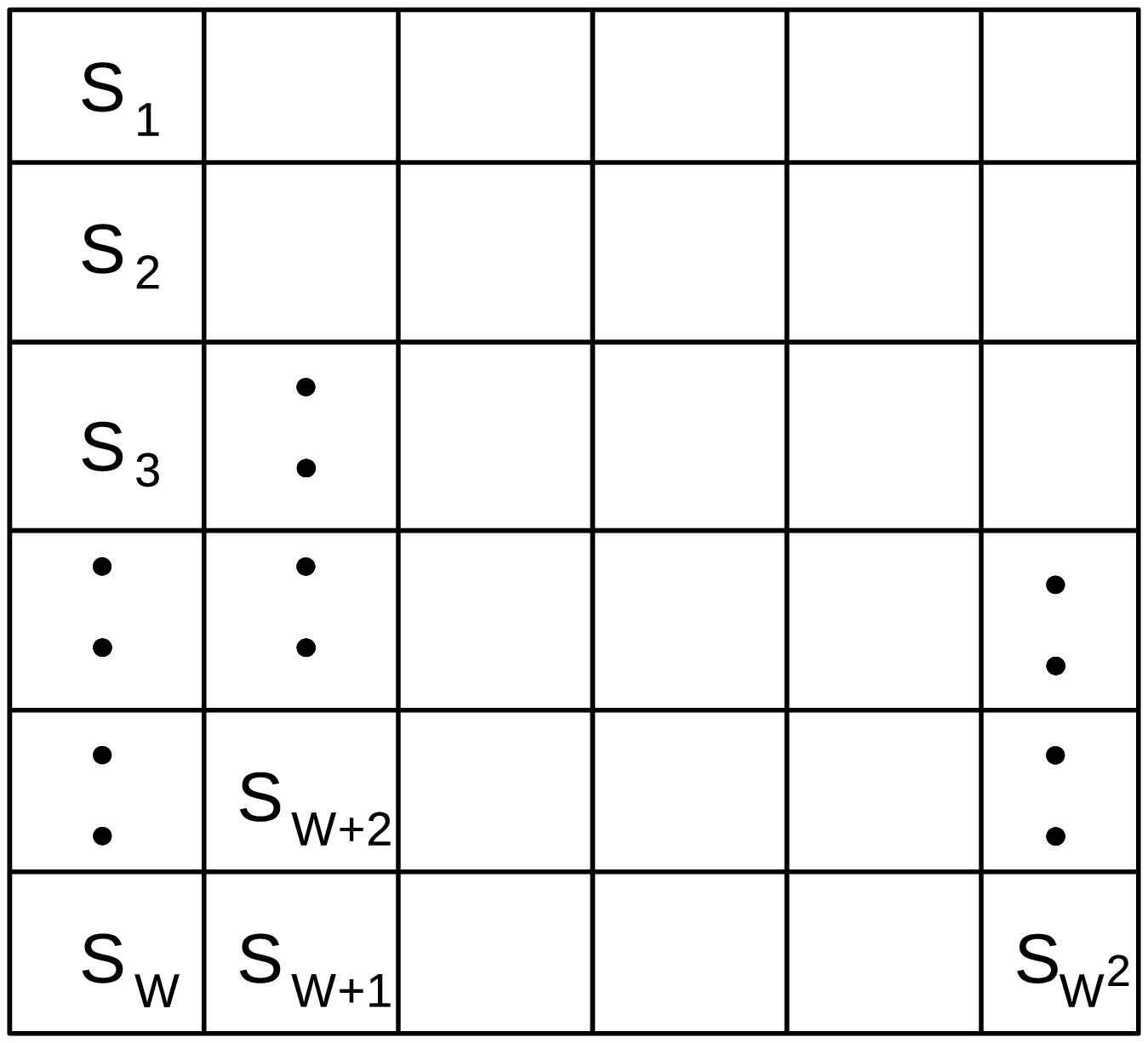}
%}
\caption{Tiling the unit square into~\(W^2 = \frac{1}{t_n^2}\) smaller~\(t_n \times t_n\) squares~\(\{S_l\}_{1 \leq l \leq W^2}.\)}
\label{fig_squares}
\end{figure}

We have a couple of preliminary definitions. Tile the unit square~\(S\) into~\(W^2\) disjoint\\\(t_n \times t_n\) squares~\(\{S_l\}_{1 \leq l \leq W^2}\)
as in Figure~\ref{fig_squares} where
\begin{equation}\label{tn_def2}
\frac{1}{20}\sqrt{\frac{M\log{n}}{n}} \leq \frac{1}{10}\sqrt{\frac{M\log(n+1)}{n+1}}  \leq \frac{r_{n+1}}{10} \leq t_n := \frac{r_{n+1}}{2\sqrt{2} + \delta_n} \leq r_{n+1}
\end{equation}
for all~\(n\) large and~\(\delta_n = [\sqrt{r_{n+1}},2\sqrt{r_{n+1}})\) is such that~\(\frac{1}{t_n}\) is an integer for all~\(n\) large. This is possible since~\(r_n \longrightarrow 0\) and so the difference~\[\frac{2\sqrt{2} + 2\sqrt{r_{n+1}}}{r_{n+1}} - \frac{2\sqrt{2} + \sqrt{r_{n+1}}}{r_{n+1}} = \frac{1}{\sqrt{r_{n+1}}} \longrightarrow \infty\] as~\(n \rightarrow \infty.\) The estimates in~(\ref{tn_def2}) follow from the fact that~\(nr_n^2 \geq M\log{n}\) for all~\(n\) large (see bounds for~\(r_n\) in the statement of the Theorem).

For~\(1 \leq j \leq n+1\) and~\(1 \leq i \leq W^2,\) let~\(E_j(i)\) be the event
that the square~\(S_i\) contains between~\(\frac{\epsilon_1 nt_n^2}{2}\) and~\(2\epsilon_2 nt_n^2\)
nodes of~\(\{X_k\}_{1 \leq k \neq j \leq n+1}\) where~\(\epsilon_1\) and~\(\epsilon_2\)
are as in~(\ref{f_eq}) and let
\begin{equation}\label{z_tot_def}
E_{dense}(n+1) := \bigcap_{j=1}^{n+1} \bigcap_{i=1}^{W^2} E_j(i).
\end{equation}
If~\(E_{dense}(n+1)\) occurs, then each~\(t_n \times t_n\)
square in~\(\{S_l\}\) contains at least one node of~\(\{X_j\}_{1 \leq j \leq n+1}.\) Since~\(2t_n \sqrt{2} < r_{n+1}\) (see~(\ref{tn_def2})),
nodes present in squares of~\(\{S_l\}\) sharing a corner are attached to each other by edges in the RGG~\(G_{n+1}.\)
Thus~\(G_{n+1}\) is connected and we let~\({\cal T}_{n+1}\) be the (unique) MST of~\(G_{n+1}.\) %as defined in~(\ref{min_weight_tree}).

For~\(1 \leq i \leq n+1,\) we denote~\(d_i\) to be the degree of the node~\(X_i\) in the MST of the component containing the node~\(X_i\)
in the graph~\(G_{n+1}\) so that if the event~\(E_{dense}(n+1)\) occurs,
then~\(d_i\) is simply the degree of the node~\(X_i\) in the MST~\({\cal T}_{n+1}.\)
Recalling that~\(\epsilon_1\) and~\(\epsilon_2\) are the bounds for the distribution~\(f\) as described in~(\ref{f_eq}), we have the following result.
\begin{Lemma}\label{diff_est} Suppose conditions~\((i)-(ii)\) in Theorem~\ref{mst_thm} hold. For all~\(n\) large,
\begin{equation}\label{e_dense_est}
\mathbb{P}(E_{dense}(n+1)) \geq 1-\frac{1}{n^3}
\end{equation}
and if~\(E_{dense}(n+1)\) occurs, then for any~\(1 \leq i \leq n+1\) the degree~\(1 \leq d_i \leq 200\epsilon_2 nr_{n+1}^2.\) Moreover, for any~\(1 \leq i \leq n+1,\) the difference
\begin{eqnarray}
|MST_{n+1}-MST(i,r_{n+1})| &\leq& \xi_{max}d_i \cdot r_{n+1}^{\alpha} \ind(E_{dense}(n+1)) \nonumber\\
&&\;\;\;\;\;\;\;\;\;+\;\;\xi_{max} n \cdot r^{\alpha}_{n+1} \ind(E^c_{dense}(n+1)),\;\;\; \nonumber\\
\label{mst_diff}
\end{eqnarray}
where~\(\xi_{max}\) is the maximum weight factor of an edge as defined prior to~(\ref{len_cyc_def}).
\end{Lemma}
%The estimate on the degree~\(d_i\) is part of the Lemma.

\emph{Proof of Lemma~\ref{diff_est}}: We first prove the estimate~(\ref{e_dense_est}).  For~\(1 \leq i \leq W^2 = \frac{1}{t_n^2},\) let~\(N_i\) be the random number of nodes of~\(\{X_k\}_{1 \leq k \leq n}\) in the~\(t_n \times t_n\) square~\(S_i.\) Since each node occurs independently with probability~\(\int_{S_i} f(x) dx\) within the square~\(S_i,\) the random variable~\(N_i\) is Binomially distributed and so from the bounds on the distribution~\(f(.)\) in~(\ref{f_eq}), the average number of nodes~\(\mathbb{E}N_i = n\int_{S_i} f(x) dx\) satisfies
\begin{equation} \label{ave_per_sq}
8 \leq \frac{\epsilon_1 M \log{n}}{20} \leq \frac{\epsilon_1nr_{n+1}^2}{10} \leq \epsilon_1 nt_n^2 \leq \mathbb{E} N_i \leq \epsilon_2 nt_n^2 \leq \epsilon_2 nr_{n+1}^2
\end{equation}
for all~\(n\) large, where~\(\epsilon_1,\epsilon_2 > 0\) are as in~(\ref{f_eq}). The estimates in~(\ref{ave_per_sq}) follow from the definition of~\(t_n\) in~(\ref{tn_def2}) and the fact that~\(nr_n^2 \geq M\log{n}\) (see statement of Theorem~\ref{mst_thm}).

Plugging the bounds obtained in~(\ref{ave_per_sq}) into the standard deviation estimates~(\ref{std_dev_up}) and~(\ref{std_dev_down}) with~\(\mu_1 = \epsilon_1 t_n^2\) and~\(\mu_2 = \epsilon_2 t_n^2\) and letting~\(\epsilon = \frac{1}{2}\) there, we therefore get
\begin{equation}
\mathbb{P}\left(\frac{\epsilon_1 nt_n^2}{2} \leq N_i \leq \frac{3\epsilon_2 nt_n^2}{2}\right) \geq 1-2\exp\left(- \frac{\epsilon_1 nt_n^2}{16}\right) \nonumber
\end{equation}
for all~\(n\) large, not depending on the choice of~\(i.\) Thus~\(\mathbb{P}(E_j(i)) \geq 1-2e^{-\frac{\epsilon_1 nt_n^2}{16}}\) and since there are~\(W^2 = \frac{1}{t_n^2}\) squares in~\(\{S_l\}\) we get from the definition of the event~\(E_{dense}(n+1)\) in~(\ref{z_tot_def}) that
\begin{eqnarray}
\mathbb{P}(E_{dense}(n+1)) &\geq& 1- (n+1)\cdot \frac{1}{t_n^2} \cdot \exp\left(-\frac{\epsilon_1 nt_n^2}{16}\right) \nonumber\\
&\geq& 1-(n+1) \cdot \frac{10}{r_{n+1}^2} \cdot \exp\left(-\frac{\epsilon_1 nr_{n+1}^2}{160}\right) \nonumber
\end{eqnarray}
since~\(t^2_n \geq \frac{r^2_{n+1}}{10}\) (see definition of~\(t_n\) in~(\ref{tn_def2})). Using the bounds in~(\ref{tn_def2}), we then get that
\begin{equation}\label{e_dense_temp}
\mathbb{P}(E_{dense}(n+1)) \geq 1- (n+1) \cdot \frac{20 n}{M\log{n}}\cdot \exp\left(-\frac{\epsilon_1 M\log{n}}{320}\right)
\end{equation}
which is at least~\(1-\frac{1}{n^3}\) for all~\(n\) large, provided~\(M > \frac{1600}{\epsilon_1}.\) We henceforth fix such an~\(M.\)

We now prove the difference estimate~(\ref{mst_diff}). Let~\({\cal T}_{n+1}\) be the weighted MST containing all the~\(n+1\) nodes~\(\{X_j\}_{1 \leq j \leq n+1}\) in the RGG~\(G_{n+1}\) with adjacency distance~\(r_{n+1}.\) Similarly, for~\(1 \leq i \leq n+1\) let~\({\cal T}(i,r_{n+1})\) be the weighted MST of the RGG~\(G(i,r_{n+1})\) formed by the nodes~\(\{X_k\}_{1 \leq k \neq i \leq n+1}\) with adjacency distance~\(r_{n+1}.\) %~\(G(i,r_{n+1}) \subseteq G_{n+1}.\)

Suppose first that the event~\(E_{dense}(n+1)\) does not occur. There are at most~\(n\) edges in either~\({\cal T}_{n+1}\) or~\({\cal T}(i,r_{n+1})\) and each such edge has length at most~\(r_{n+1}.\) Since the weight factor of any edge is at most~\(\xi_{max},\) we therefore get
\[\max(MST_{n+1}, MST(i,r_{n+1})) \ind(E^c_{dense}(n+1)) \leq \xi_{max} nr_{n+1}^{\alpha} \ind(E^c_{dense}(n+1)).\]

Henceforth we assume that~\(E_{dense}(n+1)\) occurs so that from the discussion prior to the statement of Lemma~\ref{diff_est}, we get that the RGG~\(G_{n+1}\) is connected. Arguing similarly, we get that~\(G(i,r_{n+1}) \subseteq G_{n+1}\) is also connected and wherever the location of the node~\(X_{i},\) there is always a node of~\(X_{l(i)} \in \{X_{k}\}_{1 \leq k \neq i \leq n+1}\) present within a distance of~\(t_n \sqrt{2} < r_{n+1}\) from~\(X_i.\) This means that~\({\cal T}(i,r_{n+1}) \cup \{(X_i,X_{l(i)})\}\) is a (not necessarily minimal) spanning tree containing all the nodes~\(\{X_j\}_{1 \leq j \leq n+1}.\) The edge~\((X_i,X_{l(i)})\) has weight factor at most~\(\xi_{max}\) and so the weight of the spanning tree~\({\cal T}_{n+1}\) does not exceed the weight of~\({\cal T}(i,r_{n+1})\) by more than~\(\xi_{max} r_{n+1}^{\alpha}.\) Consequently
\begin{eqnarray}
MST_{n+1} \ind(E_{dense}(n+1)) &\leq& MST(i,r_{n+1})\ind(E_{dense}(n+1))  \nonumber\\
&&\;\;\;\;\;+\;\;\xi_{max} r_{n+1}^{\alpha} \ind(E_{dense}(n+1)). \label{up1}
\end{eqnarray}

%To get an inequality in the reverse direction, we use the fact that since\\\(E_{dense}(n+1)\) occurs,
%the degree of each node in the MST~\({\cal T}_{n+1}\) is no more than a constant...~\(\frac{\pi}{6}.\) To see this is true, suppose that a node~\(u\) is attached to two nodes~\(v_1\) and~\(v_2\) in~\({\cal T}(n+1)\) such that the angle between the edges~\(e_1 = (u,v_1)\) and~\(e_2 = (u,v_2)\) is~\(\beta_{12} > \frac{\pi}{6}.\) We have that the weighted length of the edge~\(e_{12} = (v_1,v_2)\) satisfies
%\begin{eqnarray}
% wl(e_{12}) = \xi(v_1,v_2) d^{\alpha}(v_1,v_2) &=& \xi(v_1,v_2) \left(d^2(u,v_1) + d^2(u,v_2) - 2d(u,v_1)d(u,v_2)\cos{\beta_{12}}\right)^{\frac{\alpha}{2}} \nonumber\\
%&<& \xi_{max} \left(d^2(u,v_1) + d^2(u,v_2) - 2d(u,v_1)d(u,v_2)\cos{\beta_{12}}\right)^{\frac{\alpha}{2}} \\
%&=& d^(u,v_1) - d(u,v_2)(d(u,v_1) - d(u,v_2)) \nonumber
%\end{eqnarray}
%we need
%\[\xi_{max} \left(d^2(u,v_1) + d^2(u,v_2) - 2d(u,v_1)d(u,v_2)\cos{\beta_{12}}\right)^{\frac{\alpha}{2}} < \xi_{min}d^{\alpha}(u,v_1) \leq wl(e_1) \]
%which is strictly less than~\(d(u,v_1).\) This

%the Euclidean distance between~\(v_1\) CHK ABOE +etC...

\begin{figure}[tbp]
\centering
%\fbox{
\includegraphics[width=2in, trim= 60 250 150 50, clip=true]{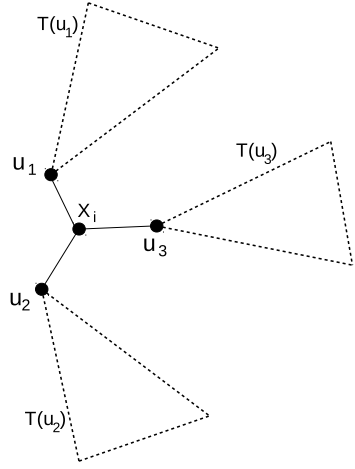}
%}
\caption{The spanning tree~\({\cal T}_{n+1}\) where the node~\(X_i\) has three neighbours~\(u_i,1 \leq i \leq 3.\) The corresponding subtrees~\({\cal T}(u_i), 1 \leq i \leq 3\) are denoted by the dotted lines.}
\label{merge_tree2}
\end{figure}

To get an inequality in the reverse direction, let~\(u_1,\ldots,u_{d_i}, d_i \geq 1\) be the neighbours of~\(X_i\) in the MST~\({\cal T}_{n+1}\) containing all the nodes~\(\{X_k\}_{1 \leq k \leq n+1}.\) Removing the node~\(X_i\) from~\({\cal T}_{n+1},\) we therefore obtain~\(d_i\) vertex disjoint trees~\({\cal T}(u_1),\ldots,{\cal T}(u_{d_i})\) with edges in~\(G(i,r_{n+1}).\)  This is illustrated in Figure~\ref{merge_tree2} where~\(d_i = 3.\)

\begin{figure}[tbp]
\centering
%\fbox{
\includegraphics[width=2in, trim= 60 250 150 50, clip=true]{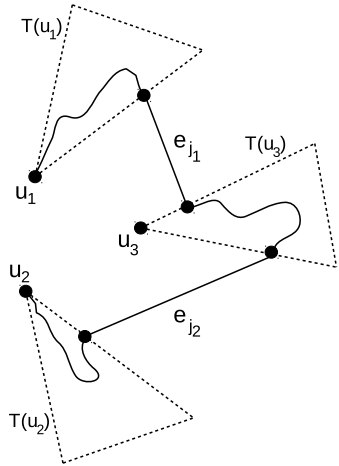}
%}
\caption{The path~\(P_{12} \subset G(i,r_{n+1})\) connecting~\(u_1\) and~\(u_2\) is formed by the union of~\(e_{j_1}, e_{j_2}\) and the wavy curves. We merge the trees~\({\cal T}(u_i), 1 \leq i \leq 3\) together by iteratively adding the bridges~\(e_{j_1}\) and~\(e_{j_2}.\)}
\label{merge_tree}
\end{figure}

Since~\(G(i,r_{n+1})\) is connected, the node~\(u_{1}\) is connected to~\(u_{2}\) by a path~\(P_{12} =(e_1,\ldots,e_g) \subseteq G(i,r_{n+1}).\) Without loss of generality, we assume that the path~\(P_{12}\) contains only edges of two types: Each edge~\(e_k\) is either an edge of a tree in~\(\{{\cal T}(u_j)\}_{1 \leq j leq d_i}\) or~\(e_k\) has one endvertex in~\({\cal T}(u_a)\) and another vertex in a different tree~\({\cal T}(u_b).\) This is illustrated in Figure~\ref{merge_tree} where~\(e_{j_1}\) is the first edge that ``leaves" the tree~\({\cal T}_{u_1}\) and therefore contains an endvertex~\(v\) not in~\({\cal T}(u_{1}).\) Since the trees~\(\{{\cal T}(u_w)\}\) contain all the nodes~\(\{X_k\}_{1 \leq k \neq i \leq n+1},\) there exists~\(1 \leq w_1 \neq u_1 \leq d_i\) such that~\(v\) belongs to~\({\cal T}(u_{w_1}).\) We denote the edge~\(e_{j_1}\) to be a bridge and adding the edge~\(e_{j_1}\) to the collection of trees~\(\{{\cal T}(u_w)\},\) we get that the graph~\(e_{j_1} \cup \cup \{{\cal T}(u_w)\}_{1 \leq w \leq d_i}\) contains~\(d_i-1\) trees.

We continue the above procedure of identifying bridges in the path~\(P_{12}\) and merging trees together until all edges of~\(P_{12}\) are exhausted. The resulting final graph~\(G^{(1)}_{merge} \subseteq G(i,r_{n+1})\) is a forest containing at most~\(d_i-1\) trees and the nodes~\(u_1\) and~\(u_2\) belong to the same tree of~\(G^{(1)}_{merge}.\)  If there still exist nodes~\(u_k\) and~\(u_l\) that belong to different trees of~\(G^{(1)}_{merge},\) we repeat the above procedure and connect~\(u_k\) and~\(u_l\) to get a new forest~\(G^{(2)}_{merge}\) containing at most~\(d_i-2\) trees. Continuing this way iteratively for a finite number of steps~\(s,\) we finally obtain a spanning tree~\(G^{(s)}_{merge}\) of~\(G(i,r_{n+1}).\)

In Figure~\ref{merge_tree}, we have illustrated the tree merging procedure described in the previous paragraph for the case when the number of neighbours~\(d_i = 3.\)  For~\(1 \leq j \leq 3,\) the dotted triangle containing the vertex~\(u_j\) represents the tree~\({\cal T}(u_j).\) The wavy curves together with the two edges~\(e_{j_1}\) and~\(e_{j_2}\) (which are bridges) form the path~\(P_{12}.\) Adding the edge~\(e_{j_1}\) merges the trees~\({\cal T}(u_1)\) and~\({\cal T}(u_3).\) Adding the edge~\(e_{j_2}\) to the resulting graph merges all the trees. In this example, all the three trees have been merged at the end of the first iteration.

The number of bridges added in the above tree merging procedure is~\(d_i-1\) (one less than the number of trees in~\(\{{\cal T}(u_w)\}_{1 \leq w \leq d_i}\)) and since any bridge has length at most~\(r_{n+1}\) and weight factor at most~\(\xi_{max},\) we get that the weight~\(W(G^{(s)}_{merge})\) of the final spanning tree of~\(G(i,r_{n+1})\) obtained satisfies
\begin{eqnarray}
MST(i,r_{n+1})\ind(E_{dense}(n+1)) &\leq& W(G^{(s)}_{merge}) \ind(E_{dense}(n+1)) \nonumber\\
&\leq& MST_{n+1}\ind(E_{dense}(n+1)) \nonumber\\
&&\;\;\;\;\;\;\;\;\;\;+\;\;\xi_{max}d_i \cdot r_{n+1}^{\alpha} \ind(E_{dense}(n+1)). \nonumber\\
\label{down1}
\end{eqnarray}
Combining~(\ref{up1}) and~(\ref{down1}) and using the fact that~\(d_i \geq 1\) if~\(E_{dense}(n+1)\) occurs (see discussion prior to statement of Lemma~\ref{diff_est}), we then get~(\ref{mst_diff}).

Finally, to upper bound~\(d_i,\) we proceed as follows. Suppose~\(X_i \in S_l\) and let~\({\cal N}_1(S_l)\) be the set of all squares in~\(\{S_k\}\) sharing a corner with~\(S_l.\) Similarly, for~\(k \geq 2\) we let~\({\cal N}_k(S_l)\) be the set of squares sharing a corner with some square in~\({\cal N}_{k-1}(S_l).\) Since every edge in~\(G_{n+1}\) has length at most~\(r_{n+1} < \frac{5t_n}{\sqrt{2}}\) (see definition of~\(t_n\) in~(\ref{tn_def2})), the neighbours of~\(X_i,\)~\(\{u_w\}_{1 \leq w \leq t},\) are all present in~\({\cal N}_4(S_l).\)

There are~\(81 < 100\) squares of~\(\{S_k\}\) in the set~\({\cal N}_{4}(S_l)\) and because~\(E_{dense}(n+1)\) occurs, every square in~\({\cal N}_4(S_l)\) has at most~\(2\epsilon_2 nt_n^2\) nodes of~\(\{X_j\}_{1 \leq j \neq i \leq n+1}\) (see definition of~\(E_{dense}(n+1)\) in~(\ref{z_tot_def})). Therefore~\(d_i \leq 200\epsilon_2 nt_n^2 \leq 200\epsilon_2 nr_{n+1}^2\) (see definition of~\(t_n\) in~(\ref{tn_def2})).~\(\qed\)

%CHK FOR BOUNDED DEGREE?? PRKVMME et+C...

%WRT MOTRE +etc...
%variance of the
%length of the minimum length cycle.\\

\subsection*{Proof of the variance estimate~(\ref{var_mst_est_main})}
We use the martingale difference method and for~\(1 \leq j \leq n+1,\) let~\({\cal F}_j = \sigma\left(\{X_k\}_{1 \leq k \leq j}\right)\) denote the~\(\sigma-\)field generated by the node positions~\(\{X_k\}_{1 \leq k \leq j}.\) Defining the martingale difference
\begin{equation}\nonumber
H_j = \mathbb{E}(MST_{n+1} | {\cal F}_j) - \mathbb{E}(MST_{n+1} | {\cal F}_{j-1}),
\end{equation}
we then have that~\(MST_{n+1} -\mathbb{E}MST_{n+1} = \sum_{j=1}^{n+1} H_j\) and so by the martingale property
\begin{equation} \nonumber
var(MST_{n+1})  = \mathbb{E}\left(\sum_{j=1}^{n+1} H_j\right)^2 = \sum_{j=1}^{n+1} \mathbb{E}H_j^2.
\end{equation}
Using the one node difference estimate from Lemma~\ref{diff_est}, we prove below that
\begin{equation}\label{gi_sec}
\sum_{j=1}^{n+1} \mathbb{E}H_j^2 \leq Cn^2r_{n+1}^{2+2\alpha}
\end{equation}
for some constant~\(C > 0\) and this obtains~(\ref{var_mst_est_main}).

\emph{Proof of~(\ref{gi_sec})}: We first rewrite the martingale difference~\(H_j\) in a more convenient form.
Let~\(\eta = (x_1,\ldots,x_{n+1})\) and~\(\eta' = (y_1,\ldots,y_{n+1})\)
be two vectors in~\((\mathbb{R}^2)^{n+1}\) so that~\(\{x_k\}_{1 \leq k \leq n+1}\)
denote the node positions in the configuration~\(\eta.\) For~\(1 \leq j \leq n,\) we let~\(\eta_j = (x_1,\ldots,x_j,y_{j+1},\ldots,y_{n+1})\) be the configuration obtained by considering the first~\(j\) entries from~\(\eta\) and the remaining from~\(\eta'.\) Let~\(G(\eta,r_{n+1})\) be the RGG formed by the nodes of~\(\eta\) with adjacency distance~\(r_{n+1}\) and
and let~\(M(\eta)\) be the total MST weight of~\(G(\eta,r_{n+1})\) as in~(\ref{min_weight_tree}). Using Fubini's theorem with the notation~\(\eta_0 = \eta'\) and~\(\eta_{n+1} = \eta,\) we then get
that the martingale difference~\[H_j = \int (M(\eta_j) - M(\eta_{j-1}))\prod_{k=j}^{n+1}f(y_k)dy_{k} \] and so we have that~\(|H_j| \leq L_j,\)
where
\begin{equation}
L_j := \int |M(\eta_j) - M(\eta_{j-1})| \prod_{k=j}^{n+1} f(y_k)dy_{k}.\label{hj_def}
\end{equation}

Recalling the event~\(E_{dense}(n+1)\) defined in~(\ref{z_tot_def})
we write~\(L_j = I_1  + I_2\) where~\(I_1\) equals
\begin{equation}
\int |M(\eta_j) - M(\eta_{j-1})|\ind(\eta_j \in E_{dense}(n+1))\ind(\eta_{j-1} \in E_{dense}(n+1)) \prod_{k=j}^{n+1} f(y_k)dy_{k} \label{i1_def_f2}
\end{equation}
and~\(I_2 = L_j-I_1\) and estimate~\(\mathbb{E}I^2_1\) and~\(\mathbb{E}I^2_2\) below, separately.

\underline{Estimate for~\(I_1\)}: Suppose the configuration~\(\eta_j \in E_{dense}(n+1)\)
and let~\(\theta_j\) be the configuration obtained by removing the node~\(x_j\) from~\(\eta_j.\)
By definition~\(\theta_j\) is also the configuration obtained by removing the node~\(y_j\) from~\(\eta_{j-1}.\)
We then recall from the discussion prior to and in Lemma~\ref{diff_est} that the RGGs formed by the nodes in~\(\eta_j\)
and in~\(\theta_j,\) with adjacency distance~\(r_{n+1},\) are both connected. Denoting the corresponding weighted MST lengths
as defined in~(\ref{min_weight_tree}) to be~\(M(\eta_j)\) and~\(M(\theta_j),\) respectively,
we get from the one node difference estimate~(\ref{mst_diff}) in Lemma~\ref{diff_est}
that
\begin{equation}\label{m_eta_j}
|M(\eta_j) - M(\theta_j)| \leq \xi_{max} d_j(\eta_j) r_{n+1}^{\alpha},
\end{equation}
where~\(1 \leq d_j(\eta_j) \leq 200 \epsilon_2 nr_{n+1}^2\) is the degree of the node~\(x_j\) in the weighted MST formed by the nodes of~\(\eta_j.\)
Similarly  we have
\begin{equation}\label{m_eta_j_p}
|M(\eta_{j-1}) - M(\theta_j)| \leq \xi_{max} d_j(\eta_{j-1}) r_{n+1}^{\alpha},
\end{equation}
where~\(d_j(\eta_{j-1})\) is the degree of the node~\(y_j\) in the weighted MST formed by the nodes of~\(\eta_{j-1}.\)

From~(\ref{m_eta_j}),~(\ref{m_eta_j_p}) and the triangle inequality, we get that~\[|M(\eta_j) - MST(\eta_{j-1})|\leq \xi_{max} (d_j(\eta_j) + d_j(\eta_{j-1}))r_{n+1}^{\alpha}.\] Using
\[(d_j(\eta_{j})+d_j(\eta_{j-1}))^2 \leq 2(d_j^2(\eta_j) + d_j^2(\eta_{j-1})) \leq 400\epsilon_2 nr_{n+1}^2 \left(d_j(\eta_j) + d_j(\eta_{j-1})\right)\]
and the definition of~\(I_1\) in~(\ref{i1_def_f2}), we then get that
\begin{eqnarray}
I^2_1 &\leq& 400\epsilon_2 \xi_{max}^2nr_{n+1}^{2+2\alpha}\int (d_j(\eta_j) + d_j(\eta_{j-1})) \prod_{k=j}^{n+1} f(y_k)dy_{k}  \nonumber\\
&=& 400\epsilon_2 \xi_{max}^2nr_{n+1}^{2+2\alpha}\left(\mathbb{E}(d_j|{\cal F}_j) + \mathbb{E}(d_j|{\cal F}_{j-1})\right) \label{i1_raw}
\end{eqnarray}
where we recall from the paragraph prior to Lemma~\ref{diff_est} that~\(d_j\) is the degree of the node~\(X_j\) in the MST of the component containing the node~\(X_j\) in the graph~\(G_{n+1}.\) Taking expectations in~(\ref{i1_raw}), we get
\begin{equation}\label{ei2}
\mathbb{E}I_1^2 \leq 800\epsilon_2 \xi_{max}^2nr_{n+1}^{2+2\alpha}\mathbb{E}(d_j).
\end{equation}

%CHNG FROM BELOW...PRKVMM +etC...

\underline{Estimate for~\(I_2\)}: The total number of edges in the spanning trees of all the components formed by the RGG containing the~\(n+1\) nodes of~\(\eta_j\) is at most~\(n\) and each such edge has length at most~\(r_{n+1}.\) Since all edge weight factors are bounded above by~\(\xi_{max},\) we use the estimate~\(\ind(A^c \cup B^c) \leq \ind(A^c) + \ind(B^c)\)
with~\[A = \{\eta_j \in E_{dense}(n+1)\} \text{ and }B = \{\eta_{j-1} \in E_{dense}(n+1)\}\]
to get that~\(I_2 \leq J_1+J_2\) where
\begin{eqnarray}
J_1 &=& nr_{n+1}^{\alpha} \xi_{max} \int \ind(\eta_j \notin E_{dense}(n+1))\prod_{k=j}^{n+1}f(y_k)dy_k \nonumber\\
&=& nr_{n+1}^{\alpha} \xi_{max} \mathbb{P}\left(E^c_{dense}(n+1)|{\cal F}_j\right) \nonumber
\end{eqnarray}
and~\(J_2 = nr_{n+1}^{\alpha} \xi_{max} \mathbb{P}\left(E^c_{dense}(n+1)|{\cal F}_{j-1}\right).\)

Thus
\[J_1^2 \leq n^2r_{n+1}^{2\alpha}\xi_{max}^2 \left(\mathbb{P}(E^c_{dense}(n+1)|{\cal F}_j)\right)^2 \leq n^2r_{n+1}^{2\alpha} \xi_{max}^2\mathbb{P}(E^c_{dense}(n+1)|{\cal F}_j)\]
and analogously, we also have~\[J_2^2 \leq n^2r_{n+1}^{2\alpha} \xi^2_{max}\mathbb{P}\left(E^c_{dense}(n+1) | {\cal F}_{j-1}\right).\]
Using~\(I_2^2 \leq (J_1+J_2)^2 \leq 2(J_1^2 + J_2^2)\) and taking expectations, we then get
\begin{equation}\label{ei22}
\mathbb{E}(I_2^2) \leq 4n^2r_{n+1}^{2\alpha}\xi^2_{max} \mathbb{P}(E^c_{dense}(n+1)) \leq 4n^2r_{n+1}^{2\alpha}\xi^2_{max}\frac{1}{n^3}
\end{equation}
using the estimate for~\(E_{dense}(n+1)\) from Lemma~\ref{diff_est}.

Adding the estimates for~\(\mathbb{E}I_1^2\) and~\(\mathbb{E}I_2^2\) in~(\ref{ei2}) and~(\ref{ei22}), respectively, we get
\[\mathbb{E}I_1^2 + \mathbb{E}I_2^2 \leq 800\epsilon_2 \xi_{max}^2nr_{n+1}^{2+2\alpha}\mathbb{E}(d_j) + 4n^2r_{n+1}^{2\alpha}\xi^2_{max}\frac{1}{n^3}\] and since~\(\mathbb{E}H_j^2 \leq \mathbb{E}L_j^2 = \mathbb{E}(I_1+I_2)^2 \leq 2(\mathbb{E}I_1^2 + \mathbb{E}I_2^2)\) we get that
\[\sum_{j=1}^{n+1} \mathbb{E}H_j^2 \leq 800\epsilon_2 \xi_{max}^2nr_{n+1}^{2+2\alpha}\mathbb{E}\left(\sum_{j=1}^{n+1}d_j\right) + 4n^2r_{n+1}^{2\alpha}\xi^2_{max}\frac{(n+1)}{n^3}.\] The sum~\(\sum_{j=1}^{n+1}d_j \leq 2n\) always since there are at most~\(n\)
edges in total among all the MSTs of the components of the RGG~\(G_{n+1}\) and so
\[\sum_{j=1}^{n+1} \mathbb{E}H_j^2 \leq \left(1600\epsilon_2 \xi_{max}^2n^2r_{n+1}^{2+2\alpha}+ 4n^2r_{n+1}^{2\alpha}\xi^2_{max}\frac{(n+1)}{n^3}\right).\]
Since~\(r_{n+1}^2 \geq \frac{M\log(n+1)}{n+1} \geq \frac{n+1}{\epsilon_2 n^3}\) for all~\(n\) large (see bounds for~\(r_n\) in the statement of the Theorem), we get that~\(\sum_{j=1}^{n+1} \mathbb{E}H_j^2 \leq 1604\epsilon_2 \xi_{max}^2n^2r_{n+1}^{2+2\alpha}\) for all~\(n\) large, proving~(\ref{gi_sec}).~\(\qed\)

\setcounter{equation}{0}
\renewcommand\theequation{\thesection.\arabic{equation}}
\section{Proof of the deviation estimates in Theorem~\ref{mst_thm}}\label{pf_mst_thm_dev}
Throughout we use Poissonization and so we construct a Poisson process~\({\cal P}\)
in the unit square~\(S\) with intensity~\(nf(.)\) as follows.
We tile the unit square~\(S\) into~\(W^2 = \frac{1}{t_n^2}\) disjoint~\(t_n \times t_n\) squares~\(\{S_i\}_{1 \leq i \leq W^2}\) (see Figure~\ref{fig_squares}), where we redefine~\(t_n = \frac{r_n}{2\sqrt{2} + \delta_n}\) and choose~\(\delta_n \in [\sqrt{r_n}, 2\sqrt{r_n})\)
such that~\(W = \frac{1}{t_n}\) is an \emph{odd} integer for all~\(n\) large. This is possible by an analogous argument following~(\ref{tn_def2}).

Each square~\(S_i, 1 \leq i \leq W^2\) has Poisson number of nodes present in it and for a formal definition we let~\(\{V_{i,k}\}_{1 \leq i \leq W^2, k\geq 1}\) be i.i.d.\ random vectors in~\(\mathbb{R}^2\) with density~\(\frac{f(x)}{\int_{S_i}f(x)dx}\ind(x \in S_i)\)
and let~\(\{N_i\}_{1 \leq i \leq W^2}\) be independent Poisson random variables
such that~\(N_i\) has mean~\(n\int_{S_i} f(x)dx\) for~\(1 \leq i \leq W^2.\)
The random variables ~\(\{N_i\}\) are independent of~\(\{V_{i,k}\}\)
and we define~\((\{V_{i,k}\}, \{N_i\})\) on the probability space~\((\Omega_0, {\cal F}_0, \mathbb{P}_0).\)

For~\(1 \leq i \leq W^2,\) we set~\(\{V_{i,k}\}_{1 \leq k \leq N_i}\) to be the nodes of~\({\cal P}\) in the square~\(S_i.\) If~\(u,v\) are any two nodes of~\({\cal P}\) that are within a Euclidean distance of~\(r_n\) of each other, we join them by an edge and denote the resulting random graph as~\(G^{(P)}_n.\) As in~(\ref{min_weight_tree}) we let~\(MST_n^{(P)}\) denote the sum of the weights of the minimum spanning trees of all components of~\(G^{(P)}_n.\) We first obtain deviation estimates for~\(MST_n^{(P)}\) and then convert via dePoissonization to obtain the corresponding estimates for~\(MST_n,\) the MST length for the Binomial process.

Analogous to~(\ref{ave_per_sq}), we have for every~\(1 \leq i \leq W^2\) that
\begin{equation} \label{ave_per_sq_poi}
8 \leq \frac{\epsilon_1 M \log{n}}{20} \leq \frac{\epsilon_1nr_{n}^2}{10} \leq \epsilon_1 nt_n^2 \leq \mathbb{E}_0 N_i \leq \epsilon_2 nt_n^2 \leq \epsilon_2 nr_n^2,
\end{equation}
where~\(\epsilon_1,\epsilon_2 > 0\) are as in~(\ref{f_eq}). If~\(E_i := \{\frac{\epsilon_1 nt_n^2}{2}  \leq N_i \leq 2\epsilon_2 nt_n^2\}\)
then using the standard deviation estimates~(\ref{std_dev_up}) and~(\ref{std_dev_down}) with~\(m = 1, \epsilon = \frac{1}{2}, \mu_1 = \epsilon_1 nt_n^2\) and~\(\epsilon_2 = nt_n^2,\) we get that~\(\mathbb{P}_0\left(E_i\right) \geq 1-\exp\left(-\frac{1}{16} \epsilon_1 nt_n^2\right) \geq 1-\exp\left(-\frac{\epsilon_1 M\log{n}}{320}\right)\)
and so if~\(E_{poi} := \bigcap_{1 \leq i \leq W^2} E_i\)
then using~\(W^2 = \frac{1}{t^2_n} \leq \frac{10}{r_n^2} \leq \frac{10n}{M\log{n}}\) (see bounds on~\(r_n\) in the statement of the theorem) we get
\begin{equation}\label{z_tot_est_p}
\mathbb{P}_0(E_{poi}) \geq 1- \frac{10M n}{\log{n}}\cdot \exp\left(-\frac{\epsilon_1 M\log{n}}{320}\right) \geq 1-\frac{1}{n^3}
\end{equation}
for all~\(n\) large, provided~\(M > \frac{1600}{\epsilon_1}.\) Henceforth, we fix such an~\(M.\)

If the event~\(E_{poi}\) occurs, then each~\(t_n \times t_n\) square in~\(\{S_l\}\)
contains a node and since~\(2t_n\sqrt{2} < r_n,\) nodes in squares of~\(\{S_l\}\) sharing a corner are connected to each other by edges
in the graph~\(G^{(P)}_n.\) Thus~\(G^{(P)}_n\) is connected.

\subsection*{Lower bounds}
We recall that we have divided the unit square~\(S\) into~\(t_n \times t_n\) squares where~\(t_n = \frac{r_n}{2\sqrt{2} + \delta_n}\) and we had chosen~\(\delta_n \in [\sqrt{r_n}, 2\sqrt{r_n})\) so that~\(W = \frac{1}{t_n}\) is an odd integer. For a real number~\(A > 0,\) we now divide each~\(t_n \times t_n\) square~\(S_i, 1 \leq i \leq W^2 = \frac{1}{t_n^2}\) into smaller disjoint~\(\frac{A(n)}{\sqrt{n}} \times \frac{A(n)}{\sqrt{n}}\) squares~\(\{R_{j}\}\) where~\(A(n) \in \left[A+ \frac{1}{(\log{n})^{1/4}}, A + \frac{2}{(\log{n})^{1/4}}\right)\) is chosen so that~\(L = \frac{t_n}{A(n)/\sqrt{n}}\) is an \emph{odd} integer for all~\(n\) large. This is possible since~\(nt_n^2 \geq \frac{\epsilon_1 M\log{n}}{20}\) (see~(\ref{ave_per_sq_poi})) and so
\begin{eqnarray}
\frac{t_n \sqrt{n}}{A + \frac{1}{(\log{n})^{1/4}}} - \frac{t_n \sqrt{n}}{A + \frac{2}{(\log{n})^{1/4}}} &=& \frac{t_n \sqrt{n}}{(\log{n})^{1/4}} \left(\frac{1}{A + \frac{1}{(\log{n})^{1/4}}}\right)\left(\frac{1}{A + \frac{2}{(\log{n})^{1/4}}}\right) \nonumber\\
&\geq& \frac{t_n \sqrt{n}}{8A^2(\log{n})^{1/4}} \nonumber
\end{eqnarray}
which~\(\longrightarrow \infty\) as~\(n \rightarrow \infty.\) For notational simplicity we henceforth refer to~\(A(n)\) simply as~\(A.\) %and treat~\(A\) as a constant.

\begin{figure}[tbp]
\centering
%\fbox{
\includegraphics[width=3in, trim= 20 180 50 110, clip=true]{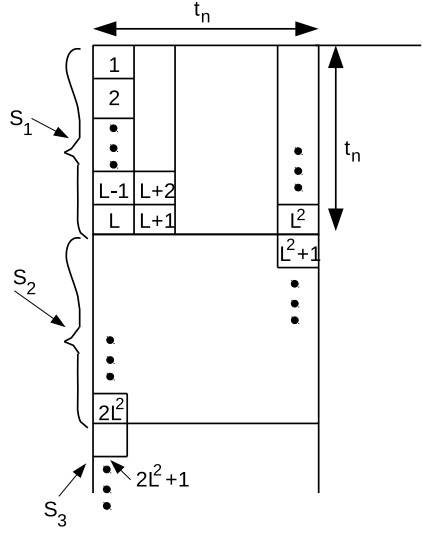}
%}
\caption{Dividing each~\(t_n \times t_n\) square into smaller subsquares.}
\label{stp_fig}
\end{figure}

The total number of~\(\frac{A}{\sqrt{n}} \times \frac{A}{\sqrt{n}}\) squares~\(\{R_j\}\) obtained as above is~\((WL)^2 = \frac{n}{A^2}\) and we label these squares in such a way that the squares~\(R_j\) and~\(R_{j+1}, 1 \leq j \leq \frac{n}{A^2}-1,\) always share a common edge. This is possible since both~\(W\) and~\(L\) are odd integers and is illustrated in Figure~\ref{stp_fig} where the square labelled~\(1\) is~\(R_1\) square labelled~\(2\) is~\(R_2\) and so on.

%As in Figure~\ref{stp_fig}, we divide each~\(t_n \times t_n\) square in~\(\{S_j\}\) into small~\(\frac{A}{\sqrt{n}} \times \frac{A}{\sqrt{n}}\) squares~\(\{R_i\}_{1 \leq i \leq \frac{n}{A^2}}\) for some fixed constant~\(A\) and assume without loss of generality that~\(\frac{\sqrt{n}}{A}\) is an integer.

The number of nodes~\(N(R_i)\) in the square~\(R_i, 1 \leq i \leq \frac{n}{A^2}\) is Poisson distributed with mean~\(\epsilon_1 A^2 \leq n\int_{R_i}f(x)dx \leq \epsilon_2 A^2,\) using~(\ref{f_eq}). Let~\(E(R_i)\) denote the event that~\(R_i\) is occupied i.e., contains at least one node of~\({\cal P},\)
and all squares sharing a corner with~\(R_i\) are empty. We remark that the event~\(E(R_i)\) is not restricted to subsquares within the same~\(t_n \times t_n\) square containing~\(R_i.\) We use the particular structure of the tiling in Figure~\ref{stp_fig} in the next subsection for the upper bound. We now recall that if the event~\(E_{poi}\) occurs, then the graph~\(G^{(P)}_n\) is connected (see discussion following~(\ref{z_tot_est_p})) and so if~\(E_{poi} \cap E(R_i)\) occurs, then there is an edge~\(e\) in the MST of~\(G^{(P)}_n\) with one endvertex in~\(R_i\) and the other endvertex in a square not sharing a corner with~\(R_i.\) The edge length of~\(e\) at least~\(\frac{A}{\sqrt{n}}\) and edge weight factor of~\(e\) is at least~\(\xi_{min},\) implying that
\begin{equation}\label{low_bound}
MST^{(P)}_n \ind(E_{poi}) \geq \frac{1}{2} \xi_{min}\sum_{i=1}^{\frac{n}{A^2}} \left(\frac{A}{\sqrt{n}}\right)^{\alpha} \ind(E(R_i) \cap E_{poi})
\end{equation}
where the factor~\(\frac{1}{2}\) occurs, since each edge is counted at most twice in the summation in~(\ref{low_bound}).

To estimate~\(H_{\alpha} := \sum_{i=1}^{\frac{n}{A^2}} \left(\frac{A}{\sqrt{n}}\right)^{\alpha} \ind(E(R_i)) ,\) we would like to split it into sums of independent r.v.s using the following construction. For a square~\(R_i,\) let~\({\cal N}(R_i)\) be the set of all squares sharing a corner with~\(R_i,\) including~\(R_i.\) If~\(R_i\) does not intersect the sides of the unit square~\(S,\) then there are~\(9\) squares in~\({\cal N}(R_i)\) and if~\(R_j\) is another square such that~\({\cal N}(R_i) \cap {\cal N}(R_j) = \emptyset,\) then the corresponding events~\(E(R_i)\) and~\(E(R_j)\) are independent, by Poisson property. We now extract nine disjoint subsets~\(\{{\cal U}_l\}_{1 \leq l \leq 9}\) of~\(\{R_i\}\) with the following properties:\\
\((A)\) If~\(R_i,R_j \in {\cal U}_l,\) then~\(\#{\cal N}(R_i) = \#{\cal N}(R_j) = 9\) and~\({\cal N}(R_i) \cap {\cal N}(R_j) = \emptyset.\)\\
\((B)\) The number of squares~\(\#{\cal U}_l \geq \frac{n}{9A^2} - \frac{4\sqrt{n}}{A}\) for each~\(1 \leq l \leq 9.\)\\
This is possible since there are at most~\(\frac{4\sqrt{n}}{A} - 4 < \frac{4\sqrt{n}}{A}\) squares in~\(\{R_k\}\) intersecting the sides of the unit square~\(S\) and the total number of squares in~\(\{R_k\}\) is~\(\frac{n}{A^2}.\)

We now write~\(H_{\alpha} = \sum_{i=1}^{\frac{n}{A^2}} \ind(E(R_i)) \geq \sum_{l=1}^{9} \sum_{R_i \in {\cal U}_l} \ind(E(R_i)),\) where each inner summation on the right side is a sum of independent Bernoulli random variables, which we bound via standard deviation estimates. Indeed for~\(1 \leq l \leq 9\) and~\(R_i \in {\cal U}_l,\) the number of nodes~\(N(R_i)\) is Poisson distributed with mean~\(n\int_{R_i}f(x)dx \in [\epsilon_1 A^2, \epsilon_2 A^2]\) (see~(\ref{f_eq})) and so~\(R_i\) is occupied with probability at least~\(1-e^{-\epsilon_1 A^2}.\) Also each of the eight squares sharing a corner with~\(R_i\) is empty with probability at least~\(e^{-\epsilon_2 A^2},\) implying that~\(\mathbb{P}(E(R_i)) \geq (1-e^{-\epsilon_1 A^2})e^{-8\epsilon_2 A^2}.\) Using the deviation estimate~(\ref{std_dev_down}) with~\(\mu_1 = (1-e^{-\epsilon_1 A^2})e^{-8\epsilon_2 A^2}, m = \frac{n}{9A^2}-\frac{4\sqrt{n}}{A}\) and~\(\epsilon = \frac{1}{m^{1/4}}\) we then get that
\begin{equation}\label{ul_est}
\mathbb{P}_0\left(\sum_{R_i \in {\cal U}_l} \ind(E(R_i)) \geq (1-\epsilon)\left(\frac{n}{9A^2}-\frac{4\sqrt{n}}{A}\right)(1-e^{-\epsilon_1 A^2})e^{-8\epsilon_2 A^2}\right) \geq 1-e^{-D_1 \epsilon^2 n}
\end{equation}
for some constant~\(D_1 > 0\) not depending on~\(l.\) Since~\(m^{1/4} < \left(\frac{n}{9A^2}\right)^{1/4},\) we get that~\(D_1\epsilon^2 n \geq 2D_2 \sqrt{n}\)  for some constant~\(D_2 > 0\) and since~\(m^{1/4} > \left(\frac{n}{10A^2}\right)^{1/4}\) for all~\(n\) large, we have \[(1-\epsilon)\left(\frac{n}{9A^2} - \frac{4\sqrt{n}}{A}\right) \geq \frac{n}{9A^2} - \frac{4\sqrt{n}}{A} - \frac{n}{A^2 m^{1/4}} \geq \frac{n}{9A^2}\left(1 - \frac{36\sqrt{A}}{n^{1/4}}\right)\] for all~\(n\) large.

Letting
\[E_{low} := \left\{H_{\alpha}\geq (1-e^{-\epsilon_1 A^2})e^{-8\epsilon_2 A^2}\frac{n}{A^2}\left(1 - \frac{36\sqrt{A}}{n^{1/4}}\right)\right\},\] we get from~(\ref{ul_est}) and a union bound that~\(\mathbb{P}_0(E_{low}) \geq 1-9e^{-2D_2 \sqrt{n}}\) and moreover, from~(\ref{low_bound}) we also get that if~\(E_{low} \cap E_{poi}\) occurs, then \[MST^{(P)}_n \geq \Delta_n := C_1(A) n^{1-\frac{\alpha}{2}} \left(1-\frac{36\sqrt{A}}{n^{1/4}}\right),\] where~\(C_1(A)\) is as defined in~(\ref{c12def}). From the estimate for the probability of the event~\(E_{low}\) above and the estimate for~\(E_{poi}\) from~(\ref{z_tot_est_p}), we therefore get~\[\mathbb{P}_0\left(MST_n^{(P)}\geq \Delta_n\right) \geq 1- \frac{1}{n^3} - 9e^{-2D_2 \sqrt{n}} \geq 1-\frac{2}{n^3}\] for all~\(n\) large.

To convert the above estimate from Poisson to the Binomial process, we let~\(E_P := \left\{MST^{(P)}_n \geq \Delta_n\right\}, E := \left\{MST_n \geq \Delta_n\right\}\) and use the dePoissonization formula
\begin{equation}\label{de_poiss_ax}
\mathbb{P}(E) \geq 1- D \sqrt{n} \mathbb{P}(E^c_P)
\end{equation}
for some constant~\(D > 0\) to get that~\(\mathbb{P}(E) \geq 1- 2D\frac{\sqrt{n}}{n^3} \geq 1-\frac{1}{n^2}\) for all~\(n\) large.
This proves the bound in~(\ref{mst_low_bounds}). Also~\[\mathbb{E}MST_n \geq \mathbb{E}MST_n \ind\left(MST_n \geq \Delta_n\right) \geq \Delta_n \left(1-\frac{1}{n^2}\right)\]
and using~\(\left(1-\frac{36\sqrt{A}}{n^{1/4}}\right)\left(1-\frac{1}{n^2}\right) \geq 1-\frac{37\sqrt{A}}{n^{1/4}}\) for all~\(n\) large, we get the lower bound on the expectation in~(\ref{exp_mst_bound}).

To prove~(\ref{de_poiss_ax}), we let~\(N_P\) denote the random number of nodes of~\({\cal P}\) in all the squares~\(\{S_j\}\) so that~\(\mathbb{E}_0 N_P = n\) and~\(\mathbb{P}_0(N_P=n) = e^{-n}\frac{n^{n}}{n!} \geq \frac{D_1}{\sqrt{n}}\) for
some constant~\(D_1 > 0,\) using the Stirling formula. Given~\(N_P = n,\) the nodes of~\({\cal P}\)
are i.i.d.\ with distribution~\(f(.)\) as defined in~(\ref{f_eq}); i.e., \(\mathbb{P}_0(E_P^c|N_P = n)  = \mathbb{P}(E^c)\) and so
\[\mathbb{P}_0(E_P^c) \geq \mathbb{P}_0(E_P^c|N_P = n) \mathbb{P}_0(N_P = n) =
\mathbb{P}(E^c) \mathbb{P}_0(N_P = n) \geq \mathbb{P}(E^c)\frac{D_1}{\sqrt{n}},\] proving~(\ref{de_poiss_ax}).~\(\qed\)

%We also need another estimate that merges trees within small squares in~\(S.\)

\subsection*{Upper bounds}
As before, we let~\({\cal P}\) be the Poisson process with intensity~\(nf(.)\) on the unit square, denote~\(G_n^{(P)}\) to be the RGG obtained from the nodes of~\({\cal P}\) with adjacency distance~\(r_n\) and set~\(MST_n^{(P)}\) to be the weighted MST length of~\(G_n^{(P)}\) as defined in~(\ref{min_weight_tree}). The main idea for obtaining the upper bound for~\(MST_n^{(P)}\) is to connect all the nodes within a square~\(R_j\) to form a subtree and then connect all these subtrees together to get an overall spanning tree. Suppose the event~\(E_{poi}\) as defined prior to~(\ref{z_tot_est_p}) occurs so that each~\(t_n \times t_n\) square in~\(\{S_l\}\) contains at least one node of the Poisson process~\({\cal P}\) and  let~\(R_{i_1},R_{i_2},\ldots,R_{i_Q}\) with~\(1 \leq i_1 < i_2 <\ldots < i_Q \leq \frac{n}{A^2}, Q \leq \frac{n}{A^2}\) be the~\(\frac{A}{\sqrt{n}} \times \frac{A}{\sqrt{n}}\) subsquares containing all the nodes of~\({\cal P}.\)

For~\(1 \leq j \leq Q,\) any two nodes in~\(R_{i_j}\) are connected by an edge (of length at most~\(\frac{A\sqrt{2}}{\sqrt{n}}\)) in the RGG~\(G_n^{(P)}\) since they belong to the same~\(t_n \times t_n\) square and~\(t_n \sqrt{2} < r_n\) (see definition of~\(t_n\) in~(\ref{tn_def2})). Therefore any spanning tree~\({\cal T}_{i_j}\) containing all the nodes of~\(R_{i_j}\) has total weight of at most~\(N(R_{i_j}) \cdot \xi_{max} \cdot \left(\frac{A\sqrt{2}}{\sqrt{n}}\right)^{\alpha},\) where~\(N(R_i)\) is the number of points of the Poisson process~\({\cal P}\) in the square~\(R_i.\) Moreover, since the event~\(E_{poi}\) occurs, the squares~\(R_{i_j}\) and~\(R_{i_{j+1}}, 1 \leq j \leq Q-1\) either belong to the same~\(t_n \times t_n\) square or belong to two different~\(t_n \times t_n\) squares sharing a common side. In any case, since~\(2t_n\sqrt{2} < r_n,\) every node in~\(R_{i_j}\) is connected to every node of~\(R_{i_{j+1}}\) by an edge of~\(G_n^{(P)}\) of length at most~\(\frac{2T_{j+1}A}{\sqrt{n}}\) where~\(T_{j+1} := i_{j+1} - i_{j}.\) We pick one such edge and call it~\(e_{j+1}.\) By definition the edge~\(e_{j+1}\) has weight at most~\(\xi_{max} \cdot \left(\frac{2T_{j+1}A}{\sqrt{n}}\right)^{\alpha}\) and for future convenience we also set~\(T_1 := i_1-1\)~and~\(T_{Q+1} := \frac{n}{A^2}-i_Q.\)

The union~\({\cal T}_{uni} := \cup_{1 \leq j \leq Q} {\cal T}_{i_j} \cup \cup_{2 \leq l \leq Q} \{e_l\}\) is a spanning tree of the graph~\(G_n\) and has weight
\begin{eqnarray}
W\left({\cal T}_{uni}\right) &\leq& \sum_{j=1}^{Q} N(R_{i_j}) \cdot \xi_{max} \cdot \left(\frac{A\sqrt{2}}{\sqrt{n}}\right)^{\alpha} +
\sum_{j=2}^{Q} \xi_{max} \cdot \left(\frac{2T_{j}A}{\sqrt{n}}\right)^{\alpha} \nonumber\\
&=& \sum_{i=1}^{\frac{n}{A^2}} N(R_{i}) \cdot \xi_{max} \cdot \left(\frac{A\sqrt{2}}{\sqrt{n}}\right)^{\alpha} +
\sum_{j=2}^{Q} \xi_{max} \cdot \left(\frac{2T_{j}A}{\sqrt{n}}\right)^{\alpha} \nonumber
\end{eqnarray}
since~\(\{R_{i_j}\}_{1 \leq j \leq Q}\) are all the squares containing at least one node of~\({\cal P}.\)
Letting\\\(Y_{\alpha} := \sum_{j=1}^{Q+1} T_j^{\alpha}\) we then get
\begin{eqnarray}\label{up_bd1}
MST_n^{(P)} \ind(E_{poi})  \leq \xi_{max} \left(\frac{2A}{\sqrt{n}}\right)^{\alpha} \left(\sum_{i=1}^{\frac{n}{A^2}} N(R_{i}) + Y_{\alpha}\right). \end{eqnarray}

The first term~\(\sum_{i=1}^{\frac{n}{A^2}} N(R_{i})\) is a Poisson random variable with mean~\(n\) since this denotes the total number of nodes of the Poisson process in the unit square. From the deviation estimate~(\ref{std_dev_up}) with~\(m=1,\mu_2 = n\) and~\(\epsilon = \frac{\log{n}}{\sqrt{n}},\) we have that~
\begin{equation}\label{e_node_est}
\mathbb{P}_0\left(\sum_{i=1}^{\frac{n}{A^2}} N(R_{i}) \leq n\left(1+\frac{\log{n}}{\sqrt{n}}\right)\right) \geq 1-e^{-C(\log{n})^2}
\end{equation}
for some constant~\(C > 0,\) using~(\ref{c_bounds}). Setting~\(E_{node} :=  \left\{\sum_{i=1}^{\frac{n}{A^2}} N(R_i) \leq n\left(1+\frac{\log{n}}{\sqrt{n}}\right)\right\},\) we get from~(\ref{up_bd1}) that
\begin{equation}\label{up_bd2}
MST_n \ind(E_{poi} \cap E_{node})  \leq \xi_{max} \left(\frac{2A}{\sqrt{n}}\right)^{\alpha} \left(n\left(1+\frac{\log{n}}{\sqrt{n}}\right)  + Y_{\alpha}\right).
\end{equation}

To evaluate the second term~\(Y_{\alpha}\) in~(\ref{up_bd2}) (which is not an i.i.d.\ sum), we consider two different cases depending on whether~\(\alpha \leq 1\) or~\(\alpha > 1.\) The term~\(Y_{\alpha}\) is well-defined for any configuration \(\omega\) of the Poisson process provided we set~\(Y_{\alpha}(\omega_0) = \left(\frac{n}{A^2}-1\right)^{\alpha}\) for the configuration~\(\omega_0\) that contains no node of the Poisson process. We now show that~\(Y_{\alpha}(\omega)\) is monotonic in~\(\omega\) in the sense that adding more points increases~\(Y_{\alpha}\) if~\(\alpha \leq 1\) and decreases~\(Y_{\alpha}\) if~\(\alpha > 1.\) This then allows us to use coupling and upper bound~\(Y_{\alpha}\) by simply considering homogenous Poisson processes.

\emph{Monotonicity of~\(Y_{\alpha}\)}: For a configuration~\(\omega \neq \omega_0\) let~\(1 \leq i_1(\omega) < \ldots < i_Q(\omega) \leq \frac{n}{A^2}\) be the indices of the squares in~\(\{R_j\}\) containing at least one node of the Poisson process~\({\cal P}.\) Letting~\(i_0(\omega) = 1\) and~\(i_{Q+1}(\omega) = \frac{n}{A^2}\) we have~\(Y_{\alpha}(\omega) = \sum_{j=0}^{Q}(i_{j+1}(\omega)-i_j(\omega))^{\alpha}.\) Suppose~\(\omega' = \omega \cup \{x\}\) is obtained by adding a single extra node at~\(x \in R_{j_0}\) for some~\(1 \leq j_0 \leq \frac{n}{A^2}.\) If~\(j_0 \in \{i_k(\omega)\}_{0 \leq k \leq Q+1},\) then~\(Y_{\alpha}(\omega') = Y_{\alpha}(\omega).\) Else there exists~\(0 \leq a \leq Q\) such that~\(i_a(\omega)  < j_0 < i_{a+1}(\omega)\) and so~\[Y_{\alpha}(\omega') = Y_{\alpha}(\omega) + (i_{a+1}(\omega)-j_0)^{\alpha} +(j_0-i_a(\omega))^{\alpha} - (i_{a+1}(\omega)-i_a(\omega))^{\alpha}.\]

If~\(\alpha \leq 1\) then using~\(a^{\alpha} + b^{\alpha} \geq (a+b)^{\alpha} \) for positive numbers~\(a,b\) we get that~\(Y_{\alpha}(\omega') \geq Y_{\alpha}(\omega).\) If~\(\alpha > 1\) then~\(a^{\alpha}  + b^{\alpha} \leq (a+b)^{\alpha}\) and so~\(Y_{\alpha}(\omega') \leq Y_{\alpha}(\omega).\) This monotonicity property together with coupling allows us to upper bound~\(Y_{\alpha}\) as follows. Letting~\(\delta = \epsilon_2\) if~\(\alpha \leq 1\) and~\(\delta = \epsilon_1\) if~\(\alpha > 1,\) we let~\({\cal P}_{\delta}\) be a homogenous Poisson process of intensity~\(\delta n\) on the unit square~\(S,\) defined on the probability space~\((\Omega_{\delta},{\cal F}_{\delta},\mathbb{P}_{\delta}).\)

Let~\(F_{\delta}\) denote the event that there is at least one node of~\({\cal P}_{\delta}\) in the unit square~\(S\) and set~\(Y^{(\delta)}_{\alpha} := \left(\frac{n}{A^2}-1\right)^{\alpha}\) if there is no node of~\({\cal P}_{\delta}\) in~\(S.\) Suppose now that~\(F_{\delta}\) occurs and as before let~\(\{i^{(\delta)}_{j}\}_{1 \leq j \leq Q_{\delta}}\) be the indices of the squares in~\(\{R_j\}\) containing at least one node of~\({\cal P}_{\delta}.\) Moreover, let~\(T^{(\delta)}_{j+1} := i^{(\delta)}_{j+1} - i^{(\delta)}_j\) for~\(1 \leq j \leq Q_{\delta}\) and set~\(T^{(\delta)}_1 := i^{(\delta)}_1-1\) and~\(T^{(\delta)}_{Q_{\delta}+1} := \frac{n}{A^2}-i^{(\delta)}_{Q_{\delta}}.\)

%for any~\(x > 0.\) %and some constant~\( D> 0\) not depending on~\(n\) or~\(x.\)

Defining~\(Y^{(\delta)}_{\alpha} := \sum_{j=1}^{Q_{\delta}+1} \left(T^{(\delta)}_j\right)^{\alpha},\) we now show for any~\(x > 0\) that
\begin{equation}\label{mon_salpha}
\mathbb{P}_{\delta}\left(Y^{(\delta)}_{\alpha} < x\right) \leq \mathbb{P}_0\left(Y_{\alpha} < x\right).
\end{equation}
\emph{Proof of~(\ref{mon_salpha})}: For~\(\alpha \leq 1\) we couple the original Poisson process~\({\cal P}\) and the homogenous process~\({\cal P}_{\delta}\) in the following way. Let~\(V_{i}, i \geq 1\) be i.i.d.\ random variables each with density~\(f(.)\) and let~\(N_V\) be a Poisson random variable with mean~\(n,\) independent of~\(\{V_i\}.\) The nodes~\(\{V_i\}_{1 \leq i \leq N_V}\) form a Poisson process with intensity~\(nf(.)\) which we denote as~\({\cal P}\) and colour green.

Let~\(U_i, i \geq 1\) be i.i.d.\ random variables each with density~\(\epsilon_2-f(.)\) where~\(\epsilon_2 \geq 1\) is as in~(\ref{f_eq}) and let~\(N_U\) be a Poisson random variable with mean~\(n(\epsilon_2-1).\) The random variables~\((\{U_i\},N_U)\) are independent of~\((\{V_i\},N_V)\) and the nodes~\(\{U_i\}_{1 \leq i  \leq N_U}\) form a Poisson process with intensity~\(n(\epsilon_2-f(.))\) which we denote as~\({\cal P}_{ext}\) and colour red. The nodes of~\({\cal P}\) and~\({\cal P}_{ext}\) together form a homogenous Poisson process with intensity~\(n\epsilon_2,\) which we denote as~\({\cal P}_{\delta}\) and define it on the probability space~\((\Omega_{\delta},{\cal F}_{\delta}, \mathbb{P}_{\delta}).\)

Let~\(\omega_{\delta} \in \Omega_{\delta}\) be any configuration and as above let~\(\{i^{(\delta)}_{j}\}_{1 \leq j \leq Q_{\delta}}\) be the indices of the squares in~\(\{R_j\}\) containing at least one node of~\({\cal P}_{\delta}\) and let~\(\{i_{j}\}_{1 \leq j \leq Q}\) be the indices of the squares in~\(\{R_j\}\) containing at least one node of~\({\cal P}.\) The indices in~\(\{i^{(\delta)}_j\}\) and~\(\{i_j\}\) depend on~\(\omega_{\delta}.\) Defining~\(Y_{\alpha} = Y_{\alpha}(\omega_{\delta})\) and~\(Y^{(\delta)}_{\alpha} = Y_{\alpha}^{(\delta)}(\omega_{\delta})\) as before, we have that~\(Y_{\alpha}\) is determined only by the green nodes of~\(\omega_{\delta}\) while~\(Y^{(\delta)}_{\alpha}\) is determined by both green and red nodes of~\(\omega_{\delta}.\)

From the monotonicity property, we therefore have that~\(Y_{\alpha}(\omega_{\delta}) \leq Y^{(\delta)}_{\alpha}(\omega_{\delta})\) and so for any~\(x > 0\) we have
\begin{equation}\label{mon_eq}
\mathbb{P}_{\delta}(Y^{(\delta)}_{\alpha} < x) \leq \mathbb{P}_{\delta}(Y_{\alpha} < x)  = \mathbb{P}_0(Y_{\alpha} < x),
\end{equation}
proving~(\ref{mon_salpha}).

If~\(\alpha > 1,\) we perform a slightly different analysis. Letting~\(\epsilon_1 \leq 1\) be as in~(\ref{f_eq}), we construct a Poisson process~\({\cal P}_{ext}\) with intensity~\(n(f(.)-\epsilon_1)\) and colour nodes of~\({\cal P}_{ext}\) red. Letting~\({\cal P}_{\delta}\) be another independent Poisson process with intensity~\(n\epsilon_1,\) we colour nodes of~\({\cal P}_{\delta}\) green. The superposition of~\({\cal P}_{ext}\) and~\({\cal P}_{\delta}\) is a Poisson process with intensity~\(nf(.),\) which we define on the probability space~\((\Omega_{\delta},{\cal F}_{\delta},\mathbb{P}_{\delta}).\) In this case, the sum~\(Y_{\alpha}\) is determined by both green and red nodes while~\(Y^{(\delta)}_{\alpha}\) is determined only by the green nodes. Again using the monotonicity property of~\(Y_{\alpha},\) we get~(\ref{mon_eq}).~\(\qed\)

To estimate~\(Y^{(\delta)}_{\alpha}\) we let~\(N^{(\delta)}(R_i),1 \leq i \leq \frac{n}{A^2},\) be the random number of nodes of~\({\cal P}_{\delta}\) in the square~\(R_i.\) The random variables~\(\{N^{(\delta)}(R_i)\}\) are i.i.d.\ Poisson distributed each with mean~\(A^2\delta.\)  For~\(i \geq \frac{n}{A^2}+1,\) we define~\(N^{(\delta)}(R_i)\) to be i.i.d.\ Poisson random variables with mean~\(A^2\delta,\) that are also independent of~\(\{N^{(\delta)}(R_i)\}_{1 \leq i \leq \frac{n}{A^2}}.\) Without loss of generality, we associate the probability measure~\(\mathbb{P}_{\delta}\) for the random variables~\(\{N^{(\delta)}(R_i)\}_{i \geq \frac{n}{A^2} +1}\) as well.

Let~\(\tilde{T}_1 := \min\{j \geq 1 : N^{(\delta)}(R_j) \geq 1\}\)  and for~\(j \geq 2,\) let~\[\tilde{T}_j := \min\{k \geq \tilde{T}_{j-1}+1 : N^{(\delta)}(R_k) \geq 1\}-\tilde{T}_{j-1}.\] The random variables~\(\{\tilde{T}_i\}\) are nearly the same as~\(\{T^{(\delta)}_i\}\) in the following sense: Suppose the event~\(F_{\delta}\) occurs so that there is at least one node of~\({\cal P}_{\delta}\) in the unit square. This means that~\(1 \leq Q_{\delta} \leq \frac{n}{A^2}\) and so~\(T^{(\delta)}_1 = i_1-1 =\tilde{T}_1-1, T^{(\delta)}_j = \tilde{T}_j\) for~\(2 \leq j \leq Q_{\delta}\) and~\(T^{(\delta)}_{Q_{\delta}+1} \leq \tilde{T}_{Q_{\delta}+1}.\) Consequently
\begin{equation}\label{dom1}
Y_{\alpha}^{(\delta)} \ind(F_{\delta}) \leq \sum_{i=1}^{Q_{\delta}+1} \tilde{T}^{\alpha}_i\ind(F_{\delta}) \leq \sum_{i=1}^{\frac{n}{A^2}+1} \tilde{T}^{\alpha}_i\ind(F_{\delta}) \leq \sum_{i=1}^{\frac{n}{A^2}+1} \tilde{T}^{\alpha}_i,
\end{equation}
since~\(Q_{\delta} \leq \frac{n}{A^2}.\)

The advantage of the above construction is that~\(\{\tilde{T}_i\}\) are i.i.d.\ geometric random variables with success parameter~\(p = 1-e^{-A^2\delta}\) and so all moments of~\(\tilde{T}_i^{\alpha}\) exist. Letting~\(\beta_i = \left(\tilde{T}_i^{\alpha} - \mathbb{E}_{\delta}\tilde{T}_i^{\alpha}\right)\) and \(\beta_{tot} = \sum_{i=1}^{\frac{n}{A^2}+1} \beta_i \) we then have that~\[\mathbb{E}_{\delta}\beta_{tot}^4  = \sum_i \mathbb{E}_{\delta} \beta_i^4 + \sum_{i < j} \mathbb{E}_{\delta} \beta_i^2 \beta_j^2  \leq C n^2\] for some constant~\(C > 0,\) by the relation~(\ref{tn_bound}) in Appendix. For~\(\epsilon > 0\) we therefore get from Chebychev's inequality that
\[\mathbb{P}_{\delta}\left(|\beta_{tot}| > \epsilon \left(\frac{n}{A^2}+1\right) \mathbb{E}_{\delta}\tilde{T}_1^{\alpha}\right) \leq D_1\frac{\mathbb{E}_{\delta}(|\beta_{tot}|^4)}{n^4\epsilon^4} \leq \frac{D_2}{n^2\epsilon^4}\]
for some constants~\(D_1,D_2 > 0.\) Setting~\(\epsilon= \frac{1}{n^{1/16}}\)we get
\begin{equation}\label{til_est1}
\mathbb{P}_{\delta}\left(\sum_{i=1}^{\frac{n}{A^2}+1} \tilde{T}^{\alpha}_i \leq  \left(1+\frac{1}{n^{1/16}}\right)\frac{n}{A^2} \mathbb{E}_{\delta}\tilde{T}_1^{\alpha}\right) \geq 1-\frac{D_2}{n^{7/4}}
\end{equation}

From the upper bound for~\(Y^{(\delta)}_{\alpha}\) in~(\ref{dom1}) and the fact that there is at least one node of~\({\cal P}_{\delta}\) in the unit square~\(S\) with probability~\(1-e^{-\delta n}\) we further get
\begin{equation}
\mathbb{P}_{\delta}\left(Y^{(\delta)}_{\alpha} \leq  \left(1+\frac{1}{n^{1/16}}\right)\frac{n}{A^2} \mathbb{E}_{\delta}\tilde{T}_1^{\alpha}\right) \geq 1-\frac{D_2}{n^{7/4}}- e^{-\delta n} \geq 1-\frac{D_3}{n^{7/4}}\nonumber
\end{equation}
for all~\(n\) large and some constant~\(D_3 > 0.\) Denoting~\(\mathbb{E}_{\delta}\tilde{T}_1^{\alpha}\) as simply~\(\mathbb{E}\tilde{T}_1^{\alpha}\) and using the coupling relation~(\ref{mon_salpha}) we finally get that
\begin{equation}
\mathbb{P}_0\left(Y_{\alpha} \leq  \left(1+\frac{1}{n^{1/16}}\right)\frac{n}{A^2} \mathbb{E}\tilde{T}_1^{\alpha}\right)\geq 1-\frac{D_3}{n^{7/4}}. \label{til_est2}
\end{equation}

Letting~\(E_{til}\) be the event on the left side of~(\ref{til_est2}) we have from~(\ref{up_bd2}) that
the term~\(MST_n^{(P)} \ind(E_{poi} \cap E_{node} \cap E_{til})  \) is bounded above by
\begin{equation}
\xi_{max} \left(\frac{2A}{\sqrt{n}}\right)^{\alpha} \left(n\left(1+\frac{\log{n}}{\sqrt{n}}\right)   + \left(1+\frac{1}{n^{1/16}}\right)\frac{n}{A^2} \mathbb{E}\tilde{T}_1^{\alpha}\right) = C_2(A) n^{1-\frac{\alpha}{2}} + bn^{-\frac{\alpha}{2}}, \label{up_bd3}
\end{equation}
where~\(C_2(A)\) is as in~(\ref{c12def}) and
\[b = \xi_{max} (2A)^{\alpha} \left(\log{n} \cdot \sqrt{n} + \left(\frac{1}{A^2}\mathbb{E}\tilde{T}_1^{\alpha}\right)n^{15/16} \right) \leq n^{16/17}\] for all~\(n\) large, again using~(\ref{tn_bound}) from Appendix. From~(\ref{up_bd3}) and the estimates for the events~\(E_{poi},E_{node}\) and~\(E_{til}\) from~(\ref{z_tot_est_p}),(\ref{e_node_est}) and~(\ref{til_est2}), respectively, we then have
\begin{equation}\label{main_poi_up}
\mathbb{P}_0\left(MST_n^{(P)} \leq C_2(A) n^{1-\frac{\alpha}{2}} + n^{\frac{16}{17}-\frac{\alpha}{2}}\right) \geq 1- \frac{1}{n^3} - e^{-D_4(\log{n})^2} - \frac{D_3}{n^{7/4}}\geq 1- \frac{D_5}{n^{7/4}}
\end{equation}
for all~\(n\) large and some constants~\(D_4,D_5 > 0.\)

\emph{Proof of~(\ref{mst_up_bounds})}: From~(\ref{main_poi_up}) and the dePoissonization formula~(\ref{de_poiss_ax}), we obtain
\begin{equation}
\mathbb{P}\left(MST_n \leq C_2(A) n^{1-\frac{\alpha}{2}} + n^{\frac{16}{17}-\frac{\alpha}{2}}\right) \geq 1- \frac{D_6 \sqrt{n}}{n^{7/4}} = 1-\frac{D_6}{n^{5/4}} \label{temp1}
\end{equation}
proving the estimate in~(\ref{mst_up_bounds}). For bounding the expectation, we let~\[\Delta_n := C_2(A) n^{1-\frac{\alpha}{2}} + n^{\frac{16}{17}-\frac{\alpha}{2}}\] and write
\begin{eqnarray}
\mathbb{E}MST_n &=& \mathbb{E}MST_n\ind(MST_n \leq \Delta_n) + \mathbb{E}MST_n \ind(MST_n > \Delta_n) \nonumber\\
&\leq& \Delta_n + \mathbb{E}MST_n \ind(MST_n > \Delta_n) \label{term2}
\end{eqnarray}
For the second term in~(\ref{term2}), we use the estimate~\(MST_n \leq n \cdot \xi_{max} \cdot r_n^{\alpha} \leq \xi_{max} n\) since each edge has length at most~\(r_n \leq 1\) (see bounds on~\(r_n\) in the statement of the Theorem) and weight at most~\(\xi_{max}.\) Using the probability estimate~(\ref{temp1}), we then get~\[\mathbb{E} MST_n \leq \Delta_n + \frac{D_6 \xi_{max} n}{n^{5/4}} \leq C_2(A)n^{1-\frac{\alpha}{2}} + 2n^{\frac{16}{17} - \frac{\alpha}{2}}\] for all~\(n\) large. This proves the expectation upper bound in~(\ref{exp_mst_bound}).~\(\qed\)

%WRT A AS C2(A) ETC... PRKVMM +ETC...

%\begin{eqnarray}
%|MST_{k+1} - MST_{k}| \leq C_1 r_{k} (\log{k}) \ind(Y_{tot}(n)) + k \sqrt{2} \ind(Y_{tot}^c(n)) \nonumber\\
%&\leq& C_2 \frac{(\log{k})^{3/2}}{\sqrt{k}} \ind(Y_{tot}(n)) + (n+1)^2 \sqrt{2} \ind(Y_{tot}^c(n)) \nonumber\\
%\label{k_dif}
%\end{eqnarray}
%for each~\(n^2 \leq k <(n+1)^2\) and for some constants~\(C_1,C_2 > 0\) not depending on~\(k\) or~\(n.\)
%The final estimate~(\ref{k_dif}) is true since from~(\ref{rn_def_mst}) we get that~\(r_k \leq C_3 \sqrt{\frac{\log{k}}{k}}\)
%for some constant~\(C_3 > 0.\)

%at the beginning of Section~\ref{pf_mst}).
%Dividing the unit square into strips of size~\(\frac{1}{\sqrt{n}} \times 1\) and arguing as in the proof of~\((b4),\)
%we obtain~\(MST_n \leq \sqrt{n} + n\frac{1}{\sqrt{n}} \sqrt{2} + 2 \leq 3\sqrt{n}\) for all~\(n\) large.

%\otimes_{k=1}^{\infty} \Omega_k

%\emph{Proof of~\((a1)-(a2)\)}:
%WRTE MORE HERE +eTC...

%We now construct an approximation~\({\cal U}_n\) of~\({\cal P}_n\) and determine the number of edges in~\({\cal U}_n.\) We write...3332

\setcounter{equation}{0}
\renewcommand\theequation{\thesection.\arabic{equation}}
\section*{Appendix}
By the continuity of the function~\(C_1(.)\) defined in~(\ref{c12def}), we have that~\(C_1(A_n) \longrightarrow C_1(A)\) as~\(n \rightarrow \infty.\)
To prove~(\ref{c_bounds}) for~\(i=2\) we let~\(T_n\) be a geometric random variable with success probability~\(p_n = 1-q_n = 1-e^{-\delta A_n^2}\) and show that
\begin{equation}\label{tn_bound}
\lim_n \mathbb{E}T_n^{\alpha} = \mathbb{E}T^{\alpha}
\end{equation}
where~\(T\) is a geometric random variable with success probability~\(p = 1-q = 1-e^{-\delta A^2}.\)
Indeed letting~\(L \geq 1\) be fixed, we write
\[\mathbb{E}T_n^{\alpha} = \sum_{k \geq 1} k^{\alpha} e^{-\delta A_n^2 (k-1)} (1-e^{-\delta A_n^2}) = Y_n(L) + Z_n(L)\]
and~\(\mathbb{E}T^{\alpha} = Y(L) + Z(L)\)
where
\begin{equation}\label{conv_one}
Y_n(L) := \sum_{1 \leq k \leq L} k^{\alpha} e^{-\delta A_n^2 (k-1)} (1-e^{-\delta A_n^2}) \longrightarrow \sum_{1 \leq k \leq L} k^{\alpha} e^{-\delta A^2 (k-1)} (1-e^{-\delta A^2}) =: Y(L)
\end{equation} as~\(n \rightarrow \infty\)
and
\[Z_n(L) := \sum_{k \geq L+1} k^{\alpha} e^{-\delta A_n^2 (k-1)} (1-e^{-\delta A_n^2}) \leq \sum_{k \geq L+1}k^{\alpha} e^{-\delta A_n^2 (k-1)}.\]
Using~\(A_n \geq A\) we further have for~\(\epsilon > 0\) that
\begin{equation}\label{conv_two}
Z_n(L) \leq \sum_{k \geq L+1}k^{\alpha} e^{-\delta A^2 (k-1)} \leq \epsilon
\end{equation}
if~\(L = L(\epsilon,\delta, \alpha,A)\) is large. For the same choice of~\(L,\) we argue as above to get that~
\begin{equation}\label{conv_three}
Z(L) := \sum_{k \geq L+1} k^{\alpha} e^{-\delta A^2 (k-1)} (1-e^{-\delta A^2}) \leq \epsilon.
\end{equation}
Finally, writing
\[|\mathbb{E}T_n^{\alpha} - \mathbb{E}T^{\alpha}| \leq |Y_n(L)-Y(L)| + |Z_n(L)- Z(L)| \leq |Y_n(L)-Y(L)| + Z_n(L)+ Z(L) \]
and using~(\ref{conv_one}),~(\ref{conv_two}) and~(\ref{conv_three}), we get for all~\(n\) large that~\(|\mathbb{E}T_n^{\alpha} - \mathbb{E}T^{\alpha}|  \leq 3\epsilon.\) Since~\(\epsilon >0\) is arbitrary, this proves~(\ref{c_bounds}) for~\(i=2.\)~\(\qed\)

\subsection*{Acknowledgement}
I thank Professors Rahul Roy, Thomas Mountford, Federico Camia, C. R. Subramanian and the referee for crucial comments that led to an improvement of the paper. I also thank IMSc for my fellowships.

\end{document}